\documentclass[12pt,a4paper]{article}

\usepackage[top=3cm,bottom=3cm,left=2.4cm,right=2.4cm]{geometry}

\usepackage[utf8]{inputenc}
\usepackage[english]{babel}
\usepackage{amsmath,amssymb,amsfonts,amsthm,graphicx}
\usepackage{hyperref}

\DeclareMathOperator*{\argmin}{arg\,min}

\begin{document}

\newcommand{\C}{\mathbb {C}}
\newcommand{\N}{\mathbb{N}}
\newcommand{\R}{\mathbb{R}}
\newcommand{\tf}{\mathcal{F}}
\newcommand{\s}{\tau}

\renewcommand{\theenumi}{\roman{enumi}}
\renewcommand{\labelenumi}{\theenumi)}

\swapnumbers

\newtheorem{0-THM1}{Theorem}[section]
\newtheorem{0-PROP1}[0-THM1]{Proposition}
\newtheorem{0-DEF1}[0-THM1]{Definition}
\newtheorem{0-PROP2}[0-THM1]{Proposition}
\newtheorem{0-EXPLE1}[0-THM1]{Example}
\newtheorem{theorem}[0-THM1]{Theorem}

\newtheorem{2-LEM1}{Lemma}[section]
\newtheorem{2-LEM2}[2-LEM1]{Lemma}
\newtheorem{2-THM1}[2-LEM1]{Theorem}
\newtheorem{2-THM2}[2-LEM1]{Theorem}

\newtheorem{3-DEF1}{Definition}[section]
\newtheorem{3-THM1}[3-DEF1]{Theorem}
\newtheorem{3-REM0}[3-DEF1]{Remark}
\newtheorem{3-DEF2}[3-DEF1]{Definition}
\newtheorem{3-LEM1}[3-DEF1]{Lemma}
\newtheorem{3-REM1}[3-DEF1]{Remark}
\newtheorem{3-COR1}[3-DEF1]{Corollary}
\newtheorem{3-REM2}[3-DEF1]{Remark}

\newtheorem{4-PROP1}{Proposition}[section]
\newtheorem{4-PROP2}[4-PROP1]{Proposition}
\newtheorem{4-REM1}[4-PROP1]{Remark}
\newtheorem{4-PROP3}[4-PROP1]{Proposition}
\newtheorem{4-PROP4}[4-PROP1]{Proposition}
\newtheorem{4-COR1}[4-PROP1]{Corollary}

\newtheorem{A-PROP1}{Proposition}[section]
\newtheorem{A-REM1}[A-PROP1]{Remark}
\newtheorem{A-PROP2}[A-PROP1]{Proposition}

\title{Stability of propagation features under time-asymptotic approximations for a class of dispersive equations}

\author{Florent Dewez\footnote{Inria, Lille-Nord Europe research center, France. E-mail: florent.dewez@inria.fr or florent.dewez@outlook.com}}

\date{}

\maketitle

\begin{abstract}
	We consider solutions in frequency bands of dispersive equations on the line defined by Fourier multipliers, these solutions being considered as wave packets. In this paper, a refinement of an existing method permitting to expand time-asymptotically the solution formulas is proposed, leading to a first term inheriting the mean position of the true solution together with a constant variance error. In particular, this first term is supported in a space-time cone whose origin position depends explicitly on the initial state, implying especially a shifted time-decay rate. This method, which takes into account both spatial and frequency information of the initial state, makes then stable some propagation features and permits a better description of the motion and the dispersion of the solutions of interest. The results are achieved firstly by making apparent the cone origin in the solution formula, secondly by applying precisely an adapted version of the stationary phase method with a new error bound, and finally by minimizing the error bound with respect to the cone origin.
\end{abstract}

\vspace{0.3cm}




\section{Introduction}

In this paper, we are interested in the time-asymptotic behaviour of wave packets of the form
\begin{equation} \label{eq:sol_formula0}
	u_f(t,x) = \frac{1}{2 \pi} \int_\R \tf u_0 (p) \, e^{-itf(p) + ixp} \, dp \; ,
\end{equation}
where $t \in \R$, $x \in \R$ and $f: \R \longrightarrow \R$ is a strictly convex symbol. We suppose that the Fourier transform $\tf u_0$ of $u_0 \in \mathcal{S}(\R)$ is supported in a bounded interval $[p_1, p_2]$, where $p_1 < p_2$ are two finite real numbers. In terms of evolution equations, wave packets of the form \eqref{eq:sol_formula0} are solutions of the following type of dispersive equations:
\begin{equation} \label{eq:evoleq0}
	\left\{ \begin{array}{l}
			\left[ i \, \partial_t - f \big(D\big) \right] u_f(t) = 0 \\ [2mm]
			u_f(0) = u_0
	\end{array} \right. \; ,
\end{equation}
for $t \in \R$, where $f(D)$ is the Fourier multiplier associated with $f$ and $u_0 \in \mathcal{S}(\R)$ the initial datum supposed to be in the frequency band $[p_1, p_2]$. For instance, the solutions of the free Schrödinger equation, of the Klein-Gordon equation or of certain higher-order evolution equations can be described by wave packets of the form \eqref{eq:sol_formula0}; we refer to \cite[Section 6]{D17-1} for further details. In the present setting, the frequency band hypothesis prevents the wave packet $u_f$ from being too much localized in space according to the uncertainty principle and makes hence challenging the task of describing its spatial propagation.\\

Some approaches solving this problem time-asymptotically have been developed. In \cite{AMHDR2012}, the authors propose to approximate the solution of the Klein-Gordon equation on a star-shaped network by a spatially localized function, the latter tending to the true solution as the time tends to infinity. This has been achieved by applying precisely the version of the stationary phase method given in \cite[Theorem 7.7.5]{H83} to an integral solution formula of the equation, the desired approximation being given by the first term of the asymptotic expansion from the stationary phase method. The principle of the stationary phase method, which consists in evaluating the integrand of the oscillatory integral of interest at the stationary point of the phase function, combined with the bounded frequency band hypothesis leads to an approximation supported in a space-time cone: this cone describes both the motion and the dispersion of the solution for large times. In particular, the results exhibit in this setting the influence of the tunnel effect on the time-decay rate of the solution.

In \cite{AMD17}, this approach has been adapted to the setting of the free Schrödinger equation on the line with initial states having integrable singular frequencies in order to study the effect of such singularities on the time-asymptotic behaviour. The version of the stationary phase method proposed in  \cite[Section 2.9]{E56} has been used since it covers the case of singular amplitudes; we mention that the authors in \cite{AMD17} propose modern formulations and detailed proofs of the results from \cite{E56}. The results show that a free particle with a singular frequency tends to travel at the speed associated with this frequency. This is highlighted by the existence of space-time cones, containing the direction given by the singular frequency, in which the time-decay rates are below the rate of the classical decay inherited from the classical dispersion.

However, the expansion provided in \cite{AMD17} is proved to blow up when approaching the space-time direction associated with the singular frequency, preventing the method from approximating uniformly the true solution in regions containing this direction. This is due to the fact that the first term of the expansion inherits the singularity of the initial state; see \cite[Sec. 3]{AMD17} for more details. To tackle this issue, another approach has been proposed in \cite{D17-1}: the precision of asymptotic expansions to one term is removed in favour of less precise but more flexible explicit and uniform estimates. In particular, they cover the above critical regions. Further this flexibility has permitted to consider not only the free Schrödinger equation but also equations of type \eqref{eq:evoleq0} with initial data having singular frequencies. The uniform estimates for the solutions have been achieved by applying a generalization \cite[Theorem 4.8]{D17-1} of the classical van der Corput Lemma \cite[Prop. 2, Chap. VIII]{S93} to the case of singular and integrable amplitudes.

We mention now that the approach developed in \cite{AMHDR2012} and the subsequent adaptations appearing in \cite{AMD17, D17-1} describe the time-asymptotic motion and the dispersion of the solutions by exploiting only the frequency information of the initial state, and not the spatial ones. For instance, this is highlighted by the fact that the origin of the space-time cones resulting from the above methods is always at the space-time point $(0, 0)$, whatever the localization of the initial state is. Consequently the associated first terms provide poor approximations of the solutions during a long time for initial states spatially far from the origin. \\

In view of this, we aim firstly at extending the approach used in \cite{AMHDR2012} to the general setting of dispersive equations of the form \eqref{eq:evoleq0} and secondly at refining it in order to exploit the spatial information of the initial datum.\\
Regarding the first point, we establish uniform and explicit remainder estimates for a time-asymptotic expansion to one term of the wave packet \eqref{eq:sol_formula0} solution of \eqref{eq:evoleq0}. As in \cite{AMHDR2012, AMD17}, we apply our new version of the stationary phase method for oscillatory integrals of the form
\begin{equation*}
	\forall \, \omega > 0 \qquad \int_\R U(p) \, e^{i \omega \psi(p)} \, dp \; ,
\end{equation*}
to the Fourier solution formula of equation \eqref{eq:evoleq0}, for a sufficiently regular and compactly supported $\tf u_0$.
In this paper, we adapt the computations of the proof for the stationary phase method in \cite{AMD17} to the case of regular and compactly supported amplitude functions $U: \R \longrightarrow \C$ and concave phase functions $\psi: \R \longrightarrow \R$. Asymptotic expansions together with uniform and explicit remainder estimates are given in Theorem \ref{2-THM1} for stationary points of the phase inside the support of the amplitude and outside the support in Theorem \ref{2-THM2}.\\
The refinement we propose is based on the two following key points:
\begin{enumerate}
	\item We establish a remainder estimate for the stationary phase method involving the $L^2$-norm of the first derivative of the amplitude, and not the $L^\infty$-norm as in the original proofs \cite{E56, AMD17}; this is done in the above mentioned Theorems \ref{2-THM1} and \ref{2-THM2}, the proof being substantially based on the application of Cauchy-Schwarz inequality to the integral representation of the remainder term.\\
	The interest of the $L^2$-norm lies in the applications to the solution formula \eqref{eq:sol_formula0}: the amplitude function $U$ being equal to the Fourier transform of the initial datum (up to a factor) in this setting, Plancherel theorem is applicable and leads to an estimate depending explicitly on the spatial part of the initial datum.
	\item We introduce an arbitrary space-time shift parametrized by a two-dimensional parameter $(t_0, x_0)$ in the integral formula \eqref{eq:sol_formula0}. Roughly speaking, this shift modifies the initial datum which is then given by the solution at time $t_0$ spatially translated by $x_0$. By applying then the above mentioned stationary phase method with the new remainder estimate, we obtain a family of time-asymptotic expansions of the solution parametrized by $(t_0, x_0)$ with explicit dependence on this parameter; in particular, the parameter $(t_0, x_0)$ is the origin of the cone in which is supported the associated first term.\\
	The parametrized family of time-asymptotic expansions for the solution of equation \eqref{eq:evoleq0} is given in Theorem \ref{3-THM1}.
\end{enumerate}
The combination of the two preceding points makes feasible the computation of the space-time parameter $(t^*, x^*)$ minimizing the remainder bound for the time-asymptotic expansion of the wave packet \eqref{eq:sol_formula0}; see Corollary \ref{3-COR1}. It is then proved in Proposition \ref{4-PROP2} that the first term associated with this optimal parameter has the same mean position as the solution; Proposition \ref{4-PROP4} shows that the difference between the variance of this approximation and the variance of the solution is an explicit constant independent from time. This refined approach permits thus to put the cone in space-time in such a way that the associated first term provides a more accurate time-asymptotic approximation of solutions in frequency bands of equations of type \eqref{eq:evoleq0}.

Let us illustrate our main result in the case of the free Schrödinger equation on the line with initial datum $u_0 \in \mathcal{S}(\R)$, namely
\begin{equation} \label{eq:free_schrodinger}
	\left\{ \begin{array}{l}
		\displaystyle i \, \partial_t u_S(t) = -\frac{1}{2} \, \partial_{xx} u_S(t) \\ [2mm]
		u_S(0) = u_0
	\end{array} \right. \; ,
\end{equation}
for $t \in \R$, whose solution is the wave function associated with the free quantum particle being in the state $u_0$ at the initial time; we note that equation \eqref{eq:free_schrodinger} is actually of the form \eqref{eq:evoleq0} with symbol $f(p) = \frac{1}{2} \, p^2$ and its solution is given by
\begin{equation} \label{eq:sol_schro}
	\forall \, (t,x) \in \R \times \R \qquad u_S(t,x) = \frac{1}{2 \pi} \int_\R \tf u_0 (p) \, e^{-\frac{1}{2} itp^2 + ixp} \, dp \; .
\end{equation}
In quantum mechanics, the frequency band hypothesis means that the particle has a momentum localized in the interval $[p_1,p_2]$. According to the physical principle of group velocity, the wave packet given by the solution will travel in space at different speeds between $p_1$ and $p_2$ over time. Hence a free wave packet in the frequency band $[p_1,p_2]$ is expected to be mainly spatially localized in an interval of the form $\big[ p_1 \, (t-t_0) + x_0 , p_2 \, (t-t_0) + x_0 \big]$, where $t_0$ and $x_0$ have to be fixed, describing hence the motion and the dispersion of the associated particle. The following result, which is a direct consequence of our main result Corollary \ref{3-COR1}, is a mathematical formulation of this principle:
\begin{0-THM1} \label{0-THM1}
	Consider the free Schrödinger equation on the line \eqref{eq:free_schrodinger} with $u_0 \in \mathcal{S}(\R)$. Let $p_1$, $p_2$, $\tilde{p}_1$ and $\tilde{p}_2$ be four finite real numbers such that $[p_1, p_2] \subset (\tilde{p}_1, \tilde{p}_2)$. Suppose $\| u_0 \|_{L^2(\R)} = 1$ and $supp \, \tf u_0 \subseteq [p_1,p_2]$, and define
	\begin{align*}
		& \bullet \quad t^* = \argmin_{\tau \in \R} \left( \int_\R x^2 \, \big| u_S(\tau, x) \big|^2 \, dx - \Big( \int_\R x \, \big| u_S(\tau, x) \big|^2 \, dx \Big)^2 \right) ; \\
		& \bullet \quad x^* = \int_\R x \, \big| u_S(t^*, x) \big|^2 \, dx \; .
	\end{align*}
	Then for all $\displaystyle (t,x) \in \left\{ (t,x) \in \big( \R \backslash \{t^* \} \big) \times \R \, \bigg| \, p_1 \leqslant \frac{x - x^*}{t - t^*} \leqslant p_2 \right\}$, we have
	\begin{align}
		& \left| u_S(t, x) - \frac{1}{\sqrt{2 \pi}} \, e^{- sgn(t-t^*) i \frac{\pi}{4}} \, e^{-it \big(\frac{x-x^*}{t - t^*}\big)^2 + ix \frac{x-x^*}{t - t^*}} \, \tf u_0 \left( \frac{x-x^*}{t - t^*} \right) |t-t^*|^{-\frac{1}{2}} \right| \nonumber \\
		 & \hspace{2cm} \leqslant C_1(\delta, \tilde{p}_1, \tilde{p}_2) \, \sqrt{\int_\R x^2 \, \big| u_S(t^*, x) \big|^2 \, dx - \Big(\int_\R x \, \big| u_S(t^*, x) \big|^2 \, dx \Big)^2} \, |t-t^*|^{-\delta} \; , \label{eq:schro_ae}
	\end{align}
	where the real number $\delta$ is arbitrarily chosen in $\big( \frac{1}{2}, \frac{3}{4} \big)$, and for all \linebreak $\displaystyle (t,x) \in \left\{ (t,x) \in \big( \R \backslash \{t^* \} \big) \times \R \, \bigg| \, \frac{x - x^*}{t - t^*} < p_1 \text{ or } p_2 < \frac{x - x^*}{t - t^*} \right\}$, we have
	\begin{align*}
		\big| u_S(t,x) \big|
			& \leqslant \Bigg( C_2(p_1, p_2, \tilde{p}_1, \tilde{p}_2) \, \sqrt{\int_\R x^2 \, \big| u_S(t^*, x) \big|^2 \, dx - \Big( \int_\R x \, \big| u_S(t^*, x) \big|^2 \, dx \Big)^2} \\
			& \hspace{1.5cm} +  C_3(p_1, p_2, \tilde{p}_1, \tilde{p}_2) \, \big\| u_0 \big\|_{L^1(\R)} \Bigg) \, |t-t^*|^{-1} \; .
	\end{align*}
	All the above constants are defined in Theorem \ref{3-THM1}.
\end{0-THM1}
\noindent See Corollary \ref{3-COR1} for the general result. Let us now make some comments on this result:
\begin{itemize}
	\item The origin of the space-time cone, in which lies the support of the first term of the expansion in \eqref{eq:schro_ae}, is actually put at the mean spatial position of the solution at the time when the variance of the solution is minimal. Hence, contrary to the preceding versions in \cite{AMHDR2012, AMD17, D17-1}, the position of the cone indeed takes into account spatial information of the solution, the cone illustrating then better the propagation and the motion of the associated particle. In particular, the mean positions of the solution and of the approximation are equal and the difference between the two variances is constant.
	\item On one hand, we observe that the first term is spatially well-localized for a solution in a narrow frequency band; on the other hand, the error is bounded by the minimal value of the standard deviation of the solution. Combined with the uncertainty principle, this exhibits a compromise: a frequency well-localized solution \eqref{eq:sol_schro} can be approximated by a function supported in a narrow space-time cone but a time sufficiently far from $t^*$ is required to achieve a good precision; on the other hand, the approximation of the solution \eqref{eq:sol_schro} with a small minimal standard deviation lies in a larger cone but the bound of the error is smaller than in the preceding case. 
	\item We remark that the time-decay rate is shifted by $t^*$, which is the time when the variance of the solution of equation \eqref{eq:free_schrodinger} is minimal; this corresponds to the fact that the origin of the cone belongs to the space-time line $\big\{ (t,x) \in \R \times \R \, \big| \, t = t^* \big\}$. Hence if we require an error smaller than a certain threshold $\varepsilon > 0$, then this precision is achieved for all $t \in \R$ satisfying
	\begin{equation*}
		|t - t^*| > C_1(\delta, \tilde{p}_1, \tilde{p}_2)^{\frac{1}{\delta}} \left( \int_\R x^2 \, \big| u_S(t^*, x) \big|^2 \, dx - \Big( \int_\R x \, \big| u_S(t^*, x) \big|^2 \, dx \Big)^2 \right)^{\frac{1}{2 \delta}} \varepsilon^{-\frac{1}{\delta}} =: \eta(\varepsilon) \; .
	\end{equation*}
	In particular if we are interested in the evolution of the solution for positive times and if $t^* < -\eta(\varepsilon)$, then the error of the approximation is smaller than $\varepsilon$ for all $t \geqslant 0$. This has to be compared with the results from the classical approach (as in \cite{AMHDR2012, AMD17}) which always imply the existence of a small time-interval with left-endpoint given by $0$ in which the error is larger than a given threshold: this is due to the lack of flexibility of the classical approach which enforces $t^* = 0$ (the decay rate is then $t^{-\frac{1}{2}}$) and puts automatically the origin of the cone at the origin of space-time.
\end{itemize}

Let us now comment on some possible improvements or applications of the present results. First of all, an interesting issue would be to apply the approach developed in this paper to more complicated settings. One may consider dispersive equations on certain networks where integral solution formulas are available, as for example the Schrödinger equation on a star-shaped network with infinite branches \cite{AMAN15} or on a tadpole graph \cite{AMAN17}. In both papers, propagation features are exhibited by exploiting wave packets in frequency bands and one may hope a better description of physical phenomena by using our refined method.

We could also consider the Schrödinger equation with a potential. In \cite{D17-2}, the time-asymptotic behaviour of the two first terms of the Dyson-Phillips series \cite[Chapter III, Theorem 1.10]{EN2000} representing the perturbed solution is studied by means of asymptotic expansions. The results concerning the second term of the series are interpreted as follows: if the initial state travels from left to right in space, then the positive frequencies of the potential tend to accelerate the motion of the second term while the negative frequencies tend to slow down or even reverse it, exhibiting advanced and retarded transmissions as well as reflections. The application of the present results could bring more information on these phenomena, in particular precise spatial information on the transmitted and reflected wave packets.

As explained in this paper, the notion of frequency band is physically meaningful and permits to describe precisely time-asymptotically the propagation of solutions of certain dispersive equations. However it is a restrictive hypothesis: for example, a function in a finite frequency band is necessarily a $\mathcal{C}^\infty$-function. Hence it would be relevant to extend this notion to functions whose Fourier transform is not necessarily compactly supported but still localized in a weaker sens. In this setting, the first term of the expansion is no longer supported in a space-time cone and so one has to quantify the localization by means of different tools. For instance, we can consider approaches based on weighted norms; such norms have been used in \cite{G07}, \cite{EHT15} or in \cite{EKMT16} to show that the continuous part of the perturbed Schrödinger evolution transports away from the origin with non-zero velocity.

Our approach makes appear naturally the shifted time-decay rate $|t-t^*|^{-\frac{1}{2}}$, where $t^*$ minimizes the variance of the solution. It would be also interesting to introduce this time-shift in other existing results to obtain greater precision. For instance, one may consider the important $L^p-L^{p'}$ estimates for which a simple argument makes apparent the shifted decay; this is proved in the following result:
\begin{0-PROP1}
	Consider the free Schrödinger equation on the line \eqref{eq:free_schrodinger} with $u_0 \in \mathcal{S}(\R)$ and define $t^* \in \R$ as follows:
	\begin{equation*}
		t^* := \argmin_{\tau \in \R} \left( \int_\R x^2 \, \big| u_S(\tau, x) \big|^2 \, dx - \Big( \int_\R x \, \big| u_S(\tau, x) \big|^2 \, dx \Big)^2 \right) .
	\end{equation*}
	Then for all $p \in [2, \infty]$, we have
	\begin{equation*}
		\forall \, t \in \R \backslash \{t^* \} \qquad \big\| u_S(t, .) \big\|_{L^p(\R)} \leqslant \left( \frac{1}{4 \pi} \right)^{-\frac{1}{2} + \frac{1}{p}} \, \big\| u_S(t^*, .) \big\|_{L^{p'}(\R)} \, |t-t^*|^{-\frac{1}{2} + \frac{1}{p}} \; ,
	\end{equation*}
	where $p'$ is the conjugate of $p$.
\end{0-PROP1}
\begin{proof}
	For the sake of clarity, we use the one-parameter group $\big(  e^{-i t \partial_{xx}} \big)_{t \in \R}$ which permits to describe the Schrödinger evolution as follows: 
	\begin{equation*}
		\forall \, t \in \R \qquad u_S(t) = e^{-i t \partial_{xx}} u_0 \; .
	\end{equation*}
	Using the group property, we have for any $t \in \R$,
	\begin{equation*}
		e^{-i t \partial_{xx}} u_0 = e^{-i (t-t^*) \partial_{xx}} \, e^{-i t^* \partial_{xx}} u_0 \; ,
	\end{equation*}
	and by applying the classical $L^p-L^{p'}$ estimate \cite[Proposition 2.2.3]{C03} to the above right-hand side, we obtain for all $t \neq t^*$,
	\begin{align*}
		\big\| u_S(t, .) \big\|_{L^p(\R)}
			& = \Big\| e^{-i t \partial_{xx}} u_0 \Big\|_{L^p(\R)} \\
			& \leqslant \left( \frac{1}{4 \pi} \right)^{-\frac{1}{2} + \frac{1}{p}} \, \Big\| e^{-i t^* \partial_{xx}} u_0 \Big\|_{L^{p'}(\R)} \, |t-t^*|^{-\frac{1}{2} + \frac{1}{p}} \\
			& = \left( \frac{1}{4 \pi} \right)^{-\frac{1}{2} + \frac{1}{p}} \, \big\| u_S(t^*, .) \big\|_{L^{p'}(\R)} \, |t-t^*|^{-\frac{1}{2} + \frac{1}{p}} \; .
	\end{align*}
	Note that we are allowed to apply the classical $L^p-L^{p'}$ estimate since $e^{-i t^* \partial_{xx}} u_0 \in \mathcal{S}(\R) \subset L^{p'}(\R)$ thanks to the hypothesis $u_0 \in \mathcal{S}(\R)$.
\end{proof}
\noindent Since the classical $L^p-L^{p'}$ estimates are exploited to establish Strichartz estimates which are themselves used to study non-linear dispersive phenomena, it is necessary to extend the above shifted $L^p-L^{p'}$ estimates to spaces larger than the Schwartz space in view of precise applications. In particular, one may examine whether $t^*$ defined above still satisfies some optimal conditions; this could be linked with the results established in \cite{CXZ2010}.

Regarding long-term perspectives of our work, one could consider the full soliton resolution for non-linear dispersive equations \cite{DKM11, DKM12-1, DKM12-2, DKM13, DKM15}, which aims at classifying the asymptotic behaviour of the non-linear solutions. A key argument for the results contained in this series of papers is the channel energy method \cite{KLLS15}, which consists in estimating the associated free solution outside a space-time cone or channel; this estimate is then used to prove that a dispersive term appearing in the decomposition of the non-linear solution goes to $0$ in the energy-space. In particular, we mention that the authors in \cite{CKS14} have to shift in time the cones and channels to derive the desired estimates. Hence one might hope that the ideas proposed in the present paper could help to understand the requirement for this shift and more generally to refine the channel energy method.

Finally we could also think about minimal escape velocities \cite{HSS00,SS87} which aim at exhi\-biting propagation features for evolution operators of type $e^{-it H}$, where $H$ is a general Hamiltonian; for instance, on may consider $H = -\partial_{xx} + V$ where $V$ is real-valued potential. As explained in \cite{HS17}, the method to establish these estimates generalizes the integration by parts which is actually crucial to describe the time-asymptotic behaviour of wave packets, as illustrated in the present paper. Our approach could bring more precision to the abstract setting and hence lead to estimates containing more information on the propagation of general wave packets.\\

The paper is organized as follows: in the following section, we begin with the new remainder estimate for an adapted version of the stationary phase method developed in \cite{E56}. We establish then time-asymptotic expansions with explicit and uniform remainder estimates depending on the shift parameter $(t_0, x_0)$ for the solution of the dispersive equation \eqref{eq:evoleq0} in Section \ref{sec:applications_de}; this section provides also the value of the optimal parameter $(t^*, x^*)$ together with the bound of the associated remainder estimate. Finally Section \ref{sec:preservation} contains results for the mean position and the variance of the first term of the time-asymptotic expansions given in Section \ref{sec:applications_de}.\\

\section{Explicit error estimates for a stationary phase me\-thod via Cauchy-Schwarz inequality} \label{sec:asympt_exp}

In this section, we establish asymptotic expansions for oscillatory integrals of the form
\begin{equation} \label{eq:oscill_int}
	\forall \, \omega > 0 \qquad \int_{\R} U(p) \, e^{i \omega \psi(p)} \, dp \; ,
\end{equation}
where the amplitude $U : \R \longrightarrow \C$ is a continuously differentiable function supported on a bounded interval and the phase $\psi : \R \longrightarrow \R$ is a strictly concave $\mathcal{C}^3$-function having a unique stationary point $p_0$. The remainder estimates we provide are explicit, uniform with respect to $p_0$ and involve the $L^2$-norm of the first derivative of the amplitude. The last point plays actually a key role in the refined method developed in Section \ref{sec:applications_de}. The asymptotic expansions together with the uniform and explicit error estimates are established in Theorems \ref{2-THM1} and \ref{2-THM2}.\\

We start by stating two technical lemmas which will be substantially used in the proof of Theorem \ref{2-THM1}.\\
The first step to expand $\omega$-asymptotically integrals of type \eqref{eq:oscill_int} consists in making simpler the phase function in order to integrate then by parts. To do so, we use the diffeomorphisms $\varphi_j$ ($j = 1,2)$ defined and studied in the following lemma. The values of these diffeomorphisms at the stationary point $p_0$ are provided in order to compute explicitly the first term of the expansions and two inequalities for $\varphi_j$ are established to estimate the errors of these expansions.\\
The proof of the following result lies mainly on an integral representation of $\varphi_j$.

\begin{2-LEM2} \label{2-LEM2}
	Let $p_0$, $\tilde{p}_1$ and $\tilde{p}_2$ be three finite real numbers such that $p_0 \in (\tilde{p}_1, \tilde{p}_2)$. Suppose that $\psi \in \mathcal{C}^3(\R, \R)$ is a strictly concave function which has a unique stationary point at $p_0$. Then, for $j = 1,2$, the function
	\begin{equation*}
		\begin{array}{ccccc}
			\varphi_j &	: &	I_j & \longrightarrow & [0, s_j] \\
			& & p & \longmapsto & \big( \psi(p_0) - \psi(p)\big)^{\frac{1}{2}} 	
		\end{array}
	\end{equation*}
	where $I_1 := [\tilde{p}_1, p_0]$, $I_2 := [p_0, \tilde{p}_2]$ and $s_j := \varphi_j(\tilde{p}_j)$, satisfies the following properties:
	\begin{enumerate}
		\item the function $\varphi_j$ is a $\mathcal{C}^2$-diffeomorphism between $I_j$ and $[0, s_j]$ ;
		\item we have
		\begin{equation*}
			\varphi_j'(p_0) = (-1)^j \sqrt{- \frac{\psi''(p_0)}{2}} \; ;
		\end{equation*}
		\item for all $p \in I_j$, the absolute value of $\varphi_j'(p)$ is lower bounded as follows:
		\begin{equation*}
			\Big| \varphi_j'(p) \Big| \geqslant \frac{1}{\sqrt{2}} \, \min_{[\tilde{p}_1, \tilde{p}_2]} \big\{-\psi''\big\} \, \big\| \psi'' \big\|_{L^{\infty}(\tilde{p}_1, \tilde{p}_2)}^{-\frac{1}{2}} \; ;
		\end{equation*}
		\item we have the following $L^{\infty}$-norm estimate for $\big( \varphi_j^{\, -1}\big)''$:
		\begin{align*}
			\Big\| \big( \varphi_j^{\, -1} \big) '' \Big\|_{L^{\infty}(0, s_j)}
			& \leqslant \big\| \psi'' \big\|_{L^{\infty}(\tilde{p}_1, \tilde{p}_2)}^{\frac{3}{2}} \, \big\| \psi^{(3)} \big\|_{L^{\infty}(\tilde{p}_1, \tilde{p}_2)} \, \min_{[\tilde{p}_1, \tilde{p}_2]} \big\{-\psi''\big\}^{-\frac{7}{2}} \\
					& \hspace{1cm} + \frac{1}{3} \, \big\| \psi'' \big\|_{L^{\infty}(\tilde{p}_1, \tilde{p}_2)}^{\frac{5}{2}} \, \big\| \psi^{(3)} \big\|_{L^{\infty}(\tilde{p}_1, \tilde{p}_2)} \, \min_{[\tilde{p}_1, \tilde{p}_2]} \big\{-\psi''\big\}^{-\frac{9}{2}} \; .
		\end{align*}
	\end{enumerate}
\end{2-LEM2}

\begin{proof}
	Let $j \in \{1,2\}$ and fix $p_0 \in (\tilde{p}_1, \tilde{p}_2)$. The proof of the present lemma is mainly based on the following integral representation of the function $\varphi_j$:
	\begin{equation} \label{eq:varphi}
		\varphi_j(p) = (-1)^j \, (p-p_0) \left( \int_0^1 \int_0^1 -\psi''\big( (1-\tau)(1-\nu)p + (\nu - \nu \tau + \tau) p_0 \big) \, (1 - \tau) \, d\nu d\tau \right)^{\frac{1}{2}} \; ,
	\end{equation}
	for all $p \in I_j$. This representation can be derived by noting firstly that
	\begin{equation*}
		\psi(p_0) - \psi(p) = \int_{p}^{p_0} \psi'(t) \, dt = -\int_{p}^{p_0} \psi'(p_0) - \psi'(t) \, dt = \int_{p}^{p_0} \int_t^{p_0} -\psi''(v) \, dv \, dt \; ;
	\end{equation*}
	then we make the change of variable $(\nu, \tau) = \big( \frac{v - t}{p_0-t}, \frac{t - p}{p_0 - p} \big)$, leading to
	\begin{equation*}
		\psi(p_0) - \psi(p) = (p - p_0)^2 \int_0^1 \int_0^1 -\psi''\big( (1-\tau)(1-\nu)p + (\nu - \nu \tau + \tau) p_0 \big) \, (1 - \tau) \, d\nu d\tau \; ,
	\end{equation*}
	and we take finally the square root of the preceding equality to obtain the desired representation \eqref{eq:varphi}.
	\begin{enumerate}
		\item Since $\psi$ is a strictly concave function on $\R$, the function $\varphi_j$ is actually the square root of the non-negative $\mathcal{C}^3$-function $p \longmapsto \psi(p_0) - \psi(p)$, showing that $\varphi_j$ is twice continuously differentiable on $I_j \backslash \{ p_0 \}$ ($\varphi_j$ is actually a $\mathcal{C}^3$-function on this domain). Let us prove that it is also twice differentiable on the whole $I_j$. To do so, note that we have for $p \in I_j \backslash \{ p_0 \}$,
		\begin{align}
			\varphi_j'(p)	& = - \frac{1}{2} \, \psi'(p) \, \big( \psi(p_0) - \psi(p) \big)^{-\frac{1}{2}} \nonumber \\
							& = -\frac{1}{2} \left( \int_p^{p_0} -\psi''(q) \, dq \right) \varphi_j(p)^{-1} \nonumber \\
							& = \frac{1}{2} \left( (p-p_0) \int_0^1 -\psi'' \big( (1-t)p + t p_0 \big) \, dt \right) \varphi_j(p)^{-1} \nonumber \\
							& = \frac{(-1)^{j}}{2} \left( \int_0^1 -\psi'' \big( (1-t)p + t p_0 \big) \, dt \right) \nonumber \\
							& \hspace{1cm} \times \left( \int_0^1 \int_0^1 -\psi''\big( (1-\tau)(1-\nu)p + (\nu - \nu \tau + \tau) p_0 \big) \, (1 - \tau) \, d\nu d\tau \right)^{-\frac{1}{2}} \; . \label{eq:deriv_varphi_2}
		\end{align}
		The preceding equality combined with the positivity of the $\mathcal{C}^1$-function $-\psi''$ shows that $\varphi_j'$ is continuously differentiable on $I_j$ whose derivative is given by
		\begin{align*}
			\varphi_j''(p)	& = \frac{(-1)^{j}}{2} \left( \int_0^1 -\psi^{(3)} \big( (p-p_0) t + p_0 \big) \, (1-t) \, dt \right) \\
							& \hspace{1cm} \times \left( \int_0^1 \int_0^1 -\psi''\big( (1-\tau)(1-\nu)p + (\nu - \nu \tau + \tau) p_0 \big) \, (1 - \tau) \, d\nu d\tau \right)^{-\frac{1}{2}} \\
							& \hspace{-0.8cm} + \frac{(-1)^{j}}{2} \left( \int_0^1 -\psi'' \big( (1-t)p + t p_0 \big) \, dt \right) \\
							& \hspace{-0.3cm} \times \left( -\frac{1}{2} \right) \frac{\int_0^1 \int_0^1 -\psi^{(3)}\big( (1-\tau)(1-\nu)p + (\nu - \nu \tau + \tau) p_0 \big) \, (1 - \tau)^2 \, (1-\nu) \, d\nu d\tau}{\left( \int_0^1 \int_0^1 -\psi''\big( (1-\tau)(1-\nu)p + (\nu - \nu \tau + \tau) p_0 \big) \, (1 - \tau) \, d\nu d\tau \right)^{\frac{3}{2}}} \; ,
		\end{align*}
		for all $p \in I_j$.\\
		Now, according to equality \eqref{eq:deriv_varphi_2}, we observe that $\varphi_j'$ is negative for $j = 1$ and positive for $j = 2$ since $-\psi'' > 0$.  By the inverse function theorem, we deduce that $\varphi_j$ is a $\mathcal{C}^2$-diffeomorphism.
		\item Thanks to the integral representation \eqref{eq:varphi}, we have
		\begin{align*}
			\varphi_j'(p_0)	& = \lim_{p \rightarrow p_0} \frac{\varphi_j(p) - \varphi_j(p_0)}{p - p_0} \\
							& = (-1)^j \, \lim_{p \rightarrow p_0} \left( \int_0^1 \int_0^1 -\psi''\big( (1-\tau)(1-\nu)p + (\nu - \nu \tau + \tau) p_0 \big) \, (1 - \tau) \, d\nu d\tau \right)^{\frac{1}{2}} \\
							& = (-1)^j \sqrt{- \frac{\psi''(p_0)}{2}} \; .
		\end{align*}
		\item From equality \eqref{eq:deriv_varphi_2} (which holds actually for all $p \in I_j$), we deduce the following lower estimate for $\varphi_j'$:
		\begin{equation} \label{eq:le_varphi_deriv}
			\forall \, p \in I_j \qquad \Big| \varphi_j'(p) \Big| \geqslant \frac{1}{\sqrt{2}} \, \min_{[\tilde{p}_1, \tilde{p}_2]} \big\{-\psi''\big\} \, \big\| \psi'' \big\|_{L^{\infty}(\tilde{p}_1, \tilde{p}_2)}^{-\frac{1}{2}} \; .
		\end{equation}
		\item From the expression of $\varphi_j''$ computed in i), we obtain the following upper estimate:
		\begin{align*}
			\forall \, p \in I_j \qquad \Big| \varphi_j''(p) \Big|
				& \leqslant \frac{1}{2 \sqrt{2}} \, \big\| \psi^{(3)} \big\|_{L^{\infty}(\tilde{p}_1, \tilde{p}_2)} \, \min_{[\tilde{p}_1, \tilde{p}_2]} \big\{-\psi''\big\}^{-\frac{1}{2}} \\
				& \hspace{1cm} + \frac{1}{6 \sqrt{2}} \, \big\| \psi'' \big\|_{L^{\infty}(\tilde{p}_1, \tilde{p}_2)} \, \big\| \psi^{(3)} \big\|_{L^{\infty}(\tilde{p}_1, \tilde{p}_2)} \, \min_{[\tilde{p}_1, \tilde{p}_2]} \big\{-\psi''\big\}^{-\frac{3}{2}} \; .
		\end{align*}
		By combining the preceding inequality with estimate \eqref{eq:le_varphi_deriv} and the following relation,
		\begin{equation*}
			\forall \, s \in [0, s_j] \qquad \big(\varphi_j^{\, -1} \big)''(s) = - \frac{\varphi_j'' \big( \varphi_j^{\, -1}(s) \big)}{\varphi_j'\big( \varphi_j^{\, -1}(s) \big)^3} \; ,
		\end{equation*}
		we obtain finally for all $s \in [0, s_j]$,
		\begin{align*}
			\Big| \big(\varphi_j^{\, -1} \big)''(s) \Big|
					& \leqslant \big\| \psi'' \big\|_{L^{\infty}(\tilde{p}_1, \tilde{p}_2)}^{\frac{3}{2}} \, \big\| \psi^{(3)} \big\|_{L^{\infty}(\tilde{p}_1, \tilde{p}_2)} \, \min_{[\tilde{p}_1, \tilde{p}_2]} \big\{-\psi''\big\}^{-\frac{7}{2}} \\
					& \hspace{1cm} + \frac{1}{3} \, \big\| \psi'' \big\|_{L^{\infty}(\tilde{p}_1, \tilde{p}_2)}^{\frac{5}{2}} \, \big\| \psi^{(3)} \big\|_{L^{\infty}(\tilde{p}_1, \tilde{p}_2)} \, \min_{[\tilde{p}_1, \tilde{p}_2]} \big\{-\psi''\big\}^{-\frac{9}{2}} \; .
		\end{align*}
	\end{enumerate}
\end{proof}

After having applied the above diffeomorphism to the integral \eqref{eq:oscill_int} (previously splitted at $p_0$), the phase becomes the simple quadratic function $s \longmapsto -s^2$. In order to make an integration by parts, creating then the first and the remainder terms of the integral, one needs an expression for a primitive of the function $s \in [0,s_0] \longmapsto e^{-i \omega s^2} \in \C$, for fixed $s_0$, $\omega > 0$. In the following lemma, a useful integral representation of such a primitive is given. As in the preceding result, its value at the origin and an inequality are also provided to compute respectively the first term of the expansion and an upper bound for the remainder term.\\
To prove Lemma \ref{2-LEM1}, we refer to the paper \cite{AMD17} which gives actually the successive primitives of more general functions by using essentially complex analysis; see \cite[Theorems 6.4, 6.5 and Corollary 6.6]{AMD17}.

\begin{2-LEM1} \label{2-LEM1}
	Let $\omega, s_0 > 0$ be two real numbers and define the function $\phi(.,\omega): [0,s_0] \longrightarrow \C$ by
	\begin{equation*}
		\phi(s,\omega) := - \int_{\Lambda(s)} e^{-i \omega z^2} \, dz \; ,
	\end{equation*}
	where $\Lambda(s)$ is the half-line in the complex plane given by
	\begin{equation*}
		\Lambda(s) := \left\{ s + t \, e^{-i \frac{\pi}{4}} \, \Big| \, t \geqslant 0 \right\} \subset \C \; .
	\end{equation*}
	Then
	\begin{enumerate}
		\item the function $\phi(.,\omega)$ is a primitive of the function $\displaystyle s \in [0,s_0] \longmapsto e^{- i \omega s^2} \in \C$ ;
		\item we have
			\begin{equation*}
				\phi(0,\omega) = - \, \frac{1}{2} \, \sqrt{\pi} \, e^{-i \frac{\pi}{4}} \, \omega^{- \frac{1}{2}} \; ;
			\end{equation*}
		\item the function $\phi(.,\omega)$ satisfies
			\begin{equation*}
				\forall \, s \in (0, s_0] \qquad \big| \phi(s,\omega) \big| \leqslant L(\delta) \, s^{1 - 2 \delta} \, \omega^{-\delta} \; ,
			\end{equation*}
			where the real number $\delta$ is arbitrarily chosen in $\big( \frac{1}{2}, 1 \big)$ and the constant $L(\delta) > 0$ is defined by
			\begin{equation*}
				L(\delta) := \frac{\sqrt{\pi}}{2} \left( \frac{1}{2 \sqrt{\pi}} \: + \: \sqrt{\frac{1}{4 \pi} + \frac{1}{2}} \right)^{2 \delta -1} \; .
			\end{equation*}
	\end{enumerate}
\end{2-LEM1}

\begin{proof}
	The function $\phi(.,\omega)$ of the present paper corresponds actually to the function $\phi_1^{(2)}(.,\omega,2,1)$ defined in \cite[Theorem 2.3]{AMD17}. Hence we apply the results established in \cite{AMD17} to the present situation:
	\begin{enumerate}
		\item One proves this first point by applying \cite[Corollary 6.6]{AMD17}, which is a consequence of Theorems 6.4 and 6.5 of \cite{AMD17}, in the case $n = 1$, $j = 2$, $\rho_j = 2$ and $\mu_j = 1$.
		\item The proof of this point lies only on basic computations which are carried out in the fourth step of the proof of \cite[Theorem 2.3]{AMD17}.
		\item The combination of Lemmas 2.4 and 2.6 of \cite{AMD17} assures this last point.
	\end{enumerate}
\end{proof}

Thanks to the two preceding lemmas, we are now in position to establish the desired asymptotic expansions with respect to the parameter $\omega$ of oscillatory integrals of type \eqref{eq:oscill_int}. In the following theorem, we are interested in the case where the stationary point $p_0$ of the phase belongs to a neighbourhood of the support of the amplitude. We emphasize that the remainder estimate we provide is different from those appearing in the original paper \cite{E56} and in \cite{AMD17}.\\
Technically speaking, we split the integral at the stationary point $p_0$ and we study separately the two resulting integrals. In each situation, the method consists firstly in using the diffeomorphism introduced in Lemma \ref{2-LEM2} to make the phase function simpler, secondly in integrating by parts to create the expansion by using Lemma \ref{2-LEM2} ii), Lemma \ref{2-LEM1} i) and ii), and finally in bounding the remainder term by combining Lemma \ref{2-LEM2} iii), iv) and Lemma \ref{2-LEM1} iii) with Cauchy-Schwarz inequality.

\begin{2-THM1} \label{2-THM1}
	Let $p_1$, $p_2$, $\tilde{p}_1$ and $\tilde{p}_2$ be four finite real numbers such that $[p_1, p_2] \subset (\tilde{p}_1, \tilde{p}_2)$. Suppose that $\psi \in \mathcal{C}^3(\R, \R): \R \longrightarrow \R$ is a strictly concave function which has a unique stationary point at $p_0 \in (\tilde{p}_1, \tilde{p}_2)$. And assume that $U \in \mathcal{C}^1(\R, \C)$ is a function satisfying
	\begin{equation*}
		supp \, U \subseteq [p_1,p_2] \; .
	\end{equation*}
	Then we have for all $\omega > 0$,
	\begin{align*}
		& \left| \int_{\R} U(p) \, e^{i \omega \psi(p)} \, dp	- \sqrt{2 \pi} \, e^{-i \frac{\pi}{4}} \, e^{i \omega \psi(p_0)} \frac{U(p_0)}{\sqrt{-\psi''(p_0)}} \, \omega^{-\frac{1}{2}} \right| \\
		& \hspace{2cm} \leqslant \Big( C_1(\psi, \delta, \tilde{p}_1, \tilde{p}_2) \, \big\| U' \big\|_{L^2(\R)} + C_2(\psi, \delta, \tilde{p}_1, \tilde{p}_2) \, \big\| U \big\|_{L^{\infty}(\R)} \Big) \, \omega^{- \delta} \; ,
	\end{align*}
	where the real number $\delta$ is arbitrarily chosen in $\big( \frac{1}{2}, \frac{3}{4} \big)$ and
	\begin{align*}
		& \bullet \quad C_1(\psi, \delta, \tilde{p}_1, \tilde{p}_2) := \frac{2^{\delta+1} \, L(\delta)}{\sqrt{3-4\delta}} \, \big( \tilde{p}_2 - \tilde{p}_1 \big)^{\frac{3 - 4 \delta}{2}} \, c_1(\psi, \delta, \tilde{p}_1, \tilde{p}_2) \; ; \\
		& \bullet \quad C_2(\psi, \delta, \tilde{p}_1, \tilde{p}_2) := \frac{2^{\delta - 1} L(\delta)}{1 - \delta} \, \big( \tilde{p}_2 - \tilde{p}_1 \big)^{2-2\delta} \, c_2(\psi, \delta, \tilde{p}_1, \tilde{p}_2) \; ; \\
		& \bullet \quad c_1(\psi, \delta, \tilde{p}_1, \tilde{p}_2) := \big\| \psi'' \big\|_{L^{\infty}(\tilde{p}_1, \tilde{p}_2)}^{\frac{3}{2}-\delta} \, \min_{[\tilde{p}_1, \tilde{p}_2]}\big\{-\psi''\big\}^{-\frac{3}{2}} \; ; \\
		&  \bullet \quad c_2(\psi, \delta, \tilde{p}_1, \tilde{p}_2) := \big\| \psi'' \big\|_{L^{\infty}(\tilde{p}_1, \tilde{p}_2)}^{\frac{5}{2}-\delta} \, \big\| \psi^{(3)} \big\|_{L^{\infty}(\tilde{p}_1, \tilde{p}_2)} \, \min_{[\tilde{p}_1, \tilde{p}_2]} \big\{-\psi''\big\}^{-\frac{7}{2}} \\
			& \hspace{4.5cm} + \frac{1}{3} \, \big\| \psi'' \big\|_{L^{\infty}(\tilde{p}_1, \tilde{p}_2)}^{\frac{7}{2}-\delta} \, \big\| \psi^{(3)} \big\|_{L^{\infty}(\tilde{p}_1, \tilde{p}_2)} \, \min_{[\tilde{p}_1, \tilde{p}_2]} \big\{-\psi''\big\}^{-\frac{9}{2}} \; .
	\end{align*}		 
	The constant $L(\delta) > 0$ is defined in Lemma \ref{2-LEM1} iii).
\end{2-THM1}

\begin{proof}
	Let $\omega > 0$ and choose $p_0 \in (\tilde{p}_1, \tilde{p}_2)$. First of all, since the support of the amplitude is included in $[p_1, p_2] \subset (\tilde{p}_1, \tilde{p}_2)$, we have clearly
	\begin{equation*}
		\int_\R U(p) \, e^{i \omega \psi(p)} \, dp \, = \int_{\tilde{p}_1}^{\tilde{p}_2} U(p) \, e^{i \omega \psi(p)} \, dp \, =: I(\omega) \; .
	\end{equation*}
	Splitting the above integral at the point $p_0$ and using the two $\mathcal{C}^2$-diffeomorphisms defined in Lemma \ref{2-LEM2}, we obtain
	\begin{align*}
		I(\omega)	& = - \int_0^{s_1} \big(U \circ \varphi_1^{\, -1} \big)(p) \, \big(\varphi_1^{\, -1}\big)'(p) \, e^{-i \omega s^2} \, ds \, e^{i \omega \psi(p_0)} \\
					& \hspace{1cm} + \int_0^{s_2} \big(U \circ \varphi_2^{\, -1} \big)(p) \, \big(\varphi_2^{\, -1}\big)'(p) \, e^{-i \omega s^2} \, ds \, e^{i \omega \psi(p_0)} \; ;
	\end{align*}
	note that we have used the fact that $\varphi_1$ and $\varphi_2$ are respectively decreasing and increasing. We integrate now by parts by using the primitive $s \longmapsto \phi(s,\omega)$ given in Lemma \ref{2-LEM1} and the regularity of $\varphi_j$:
	\begin{align*}
		& (-1)^j \int_0^{s_j} \big(U \circ \varphi_j^{\, -1} \big)(s) \, \big(\varphi_j^{\, -1}\big)'(s) \, e^{-i \omega s^2} \, ds \\
		& \hspace{1.5cm} = (-1)^j \Big[ \big(U \circ \varphi_j^{\, -1} \big)(s) \, \big(\varphi_j^{\, -1}\big)'(s) \, \phi(s,\omega) \Big]_0^{s_j} \\
		& \hspace{3cm} + (-1)^{j+1} \int_0^{s_j} \Big( \big(U \circ \varphi_j^{\, -1} \big) \, \big(\varphi_j^{\, -1}\big)'\Big)'(s) \, \phi(s, \omega) \, ds \\
		& \hspace{1.5cm} = (-1)^{j+1} \big(U \circ \varphi_j^{\, -1} \big)(0) \, \big(\varphi_j^{\, -1}\big)'(0) \, \phi(0,\omega) \\
		& \hspace{3cm} + (-1)^{j+1} \int_0^{s_j} \Big( \big(U \circ \varphi_j^{\, -1} \big) \, \big(\varphi_j^{\, -1}\big)'\Big)'(s) \, \phi(s, \omega) \, ds \\
		& \hspace{1.5cm} = \frac{1}{2} \sqrt{2 \pi} \, e^{-i \frac{\pi}{4}} \frac{U(p_0)}{\sqrt{- \psi''(p_0)}} \, \omega^{-\frac{1}{2}} \\
		& \hspace{3cm} + (-1)^{j+1} \int_0^{s_j} \Big( \big(U \circ \varphi_j^{\, -1} \big) \, \big(\varphi_j^{\, -1}\big)'\Big)'(s) \, \phi(s, \omega) \, ds \; ;
	\end{align*}
	the second equality has been obtained by using the fact that $U(\tilde{p}_j) = 0$ and the last one by applying Lemma \ref{2-LEM2} ii) and Lemma \ref{2-LEM1} ii). Hence it follows
	\begin{align*}
		I(\omega)	& = \sqrt{2 \pi} \, e^{-i \frac{\pi}{4}} \, e^{i \omega \psi(p_0)} \, \frac{U(p_0)}{\sqrt{- \psi''(p_0)}} \, \omega^{-\frac{1}{2}} \\
					& \hspace{1.5cm} + \sum_{j=1}^2 (-1)^{j+1} \int_0^{s_j} \Big( \big(U \circ \varphi_j^{\, -1} \big) \, \big(\varphi_j^{\, -1}\big)'\Big)'(s) \, \phi(s, \omega) \, ds \, e^{i \omega \psi(p_0)} \; .
	\end{align*}
	To estimate each term of the remainder, we proceed as follows:
	\begin{align*}
		& \left| (-1)^{j+1} \int_0^{s_j} \Big( \big(U \circ \varphi_j^{\, -1} \big) \, \big(\varphi_j^{\, -1}\big)'\Big)'(s) \, \phi(s, \omega) \, ds \right| \\
		& \hspace{1.5cm} \leqslant \left| \int_0^{s_j} \big(U' \circ \varphi_j^{\, -1} \big)(s) \, \big(\varphi_j^{\, -1}\big)'(s)^2 \, \phi(s, \omega) \, ds \right| \\
		& \hspace{3cm} + \left| \int_0^{s_j} \big(U \circ \varphi_j^{\, -1} \big)(s) \, \big(\varphi_j^{\, -1}\big)''(s) \, \phi(s, \omega) \, ds \right| \\
		& \hspace{1.5cm} \leqslant \left( \int_0^{s_j} \Big| \big(U' \circ \varphi_j^{\, -1} \big)(s) \, \big(\varphi_j^{\, -1}\big)'(s)^2 \Big|^2 \, ds \right)^{\frac{1}{2}} \left( \int_0^{s_j} \big| \phi(s, \omega) \big|^2 \, ds \right)^{\frac{1}{2}}  \\
		& \hspace{3cm} + \int_0^{s_j} \big| \phi(s, \omega) \big|\,  ds \, \big\| U \big\|_{L^{\infty}(\R)} \, \Big\| \big(\varphi_j^{\, -1}\big)'' \Big\|_{L^{\infty}(0,s_j)} \; ;
	\end{align*}
	let us remark that we have applied Cauchy-Schwarz inequality to the first integral. We continue the proof by estimating each resulting term; first of all, by making the change of variable $p = \varphi_j^{\, -1}(s)$ and by using Lemma \ref{2-LEM2} iii), we obtain
	\begin{equation*}
		 \left( \int_0^{s_j} \Big| \big(U' \circ \varphi_j^{\, -1} \big)(s) \, \big(\varphi_j^{\, -1}\big)'(s)^2 \Big|^2 \, ds \right)^{\frac{1}{2}} \leqslant 2^{\frac{3}{4}} \, \big\| \psi'' \big\|_{L^{\infty}(\tilde{p}_1, \tilde{p}_2)}^{\frac{3}{4}} \, \min_{[\tilde{p}_1, \tilde{p}_2]}\big\{-\psi''\big\}^{-\frac{3}{2}} \, \big\| U' \big\|_{L^2(\R)} \; .
	\end{equation*}
	Then we use the point iii) of Lemma \ref{2-LEM1} to derive the two following inequalities:
	\begin{align*}
		& \bullet \quad \int_0^{s_j} \big| \phi(s, \omega) \big|\, ds \leqslant L(\delta) \int_0^{s_j} s^{1-2\delta} \, ds \, \omega^{-\delta} \leqslant \frac{L(\delta)}{2 - 2 \delta} \, \varphi_j(\tilde{p}_j)^{2-2\delta} \, \omega^{-\delta} \; ;\\
		& \bullet \quad \left( \int_0^{s_j} \big| \phi(s, \omega) \big|^2 \, ds \right)^{\frac{1}{2}} \leqslant \frac{L(\delta)}{\sqrt{3-4\delta}} \, \varphi_j(\tilde{p}_j)^{\frac{3-4\delta}{2}} \, \omega^{-\delta} \; .
	\end{align*}
	By using the integral representation \eqref{eq:varphi} of $\varphi_j$, we obtain
	\begin{equation*}
		\varphi_j(\tilde{p}_j) \leqslant \frac{1}{\sqrt{2}} \, \big\| \psi'' \big\|_{L^{\infty}(\tilde{p}_1, \tilde{p}_2)}^{\frac{1}{2}} \, (\tilde{p_2} - \tilde{p}_1) \; ,
	\end{equation*}
	which permits to deduce
	\begin{align*}
		& \bullet \quad \int_0^{s_j} \big| \phi(s, \omega) \big|\, ds \leqslant \frac{1}{2^{1-\delta}} \, \frac{L(\delta)}{2 - 2 \delta} \big( \tilde{p}_2 - \tilde{p}_1 \big)^{2-2\delta} \, \big\| \psi'' \big\|_{L^{\infty}(\tilde{p}_1, \tilde{p}_2)}^{1 - \delta} \, \omega^{-\delta} \; ;\\
		& \bullet \quad \left( \int_0^{s_j} \big| \phi(s, \omega) \big|^2 \, ds \right)^{\frac{1}{2}} \leqslant \frac{1}{2^{\frac{3}{4} - \delta}} \, \frac{L(\delta)}{\sqrt{3-4\delta}} \, \big( \tilde{p}_2 - \tilde{p}_1 \big)^{\frac{3-4\delta}{2}} \, \big\| \psi'' \big\|_{L^{\infty}(\tilde{p}_1, \tilde{p}_2)}^{\frac{3}{4}-\delta} \, \omega^{-\delta} \; .
	\end{align*}
	And, from Lemma \ref{2-LEM2} iv), we recall that
	\begin{align*}
		 \Big\| \big( \varphi_j^{\, -1} \big)'' \Big\|_{L^{\infty}(0, s_j)}
		 	& \leqslant \big\| \psi'' \big\|_{L^{\infty}(\tilde{p}_1, \tilde{p}_2)}^{\frac{3}{2}} \, \big\| \psi^{(3)} \big\|_{L^{\infty}(\tilde{p}_1, \tilde{p}_2)} \, \min_{[\tilde{p}_1, \tilde{p}_2]} \big\{-\psi''\big\}^{-\frac{7}{2}} \\
			& \hspace{1cm} + \frac{1}{3} \, \big\| \psi'' \big\|_{L^{\infty}(\tilde{p}_1, \tilde{p}_2)}^{\frac{5}{2}} \, \big\| \psi^{(3)} \big\|_{L^{\infty}(\tilde{p}_1, \tilde{p}_2)} \, \min_{[\tilde{p}_1, \tilde{p}_2]} \big\{-\psi''\big\}^{-\frac{9}{2}} \; .
	\end{align*}
	Putting everything together provides the desired estimate, namely,
	\begin{align*}
		& \left| I(\omega) - \sqrt{2 \pi} \, e^{-i \frac{\pi}{4}} \, e^{i \omega \psi(p_0)} \, \frac{U(p_0)}{\sqrt{- \psi''(p_0)}} \, \omega^{-\frac{1}{2}} \right| \\
		& \hspace{1cm} \leqslant \sum_{j=1}^2 \left| (-1)^{j+1} \int_0^{s_j} \Big( \big(U \circ \varphi_j^{\, -1} \big) \, \big(\varphi_j^{\, -1}\big)'\Big)'(s) \, \phi(s, \omega) \, ds \, e^{i \omega \psi(p_0)} \right| \\
		& \hspace{1cm} \leqslant \frac{2^{\delta+1} \, L(\delta)}{\sqrt{3-4\delta}} \, \big( \tilde{p}_2 - \tilde{p}_1 \big)^{\frac{3 - 4 \delta}{2}} \, \big\| \psi'' \big\|_{L^{\infty}(\tilde{p}_1, \tilde{p}_2)}^{\frac{3}{2}-\delta} \, \min_{[\tilde{p}_1, \tilde{p}_2]} \big\{-\psi''\big\}^{-\frac{3}{2}} \, \big\| U' \big\|_{L^2(\R)} \, \omega^{-\delta} \\
		& \hspace{2cm} + \frac{2^{\delta - 1} L(\delta)}{1 - \delta} \, \big( \tilde{p}_2 - \tilde{p}_1 \big)^{2-2\delta} \bigg( \big\| \psi'' \big\|_{L^{\infty}(\tilde{p}_1, \tilde{p}_2)}^{\frac{5}{2}-\delta} \, \big\| \psi^{(3)} \big\|_{L^{\infty}(\tilde{p}_1, \tilde{p}_2)} \, \min_{[\tilde{p}_1, \tilde{p}_2]} \big\{-\psi''\big\}^{-\frac{7}{2}} \\
		& \hspace{3cm} + \frac{1}{3} \, \big\| \psi'' \big\|_{L^{\infty}(\tilde{p}_1, \tilde{p}_2)}^{\frac{7}{2}-\delta} \, \big\| \psi^{(3)} \big\|_{L^{\infty}(\tilde{p}_1, \tilde{p}_2)} \, \min_{[\tilde{p}_1, \tilde{p}_2]} \big\{-\psi''\big\}^{-\frac{9}{2}} \bigg) \, \| U \|_{L^{\infty}(\R)} \, \omega^{-\delta} \; .
	\end{align*}
\end{proof}

We end this section by providing an explicit and uniform bound for oscillatory integrals of type \eqref{eq:oscill_int} in the case where there is no stationary point inside the support of the amplitude, making the decay with respect to $\omega$ faster. As above, the estimate involves the $L^2$-norm of the first derivative of the amplitude in view of applications to dispersive equations in the following section.\\
The proof of the following result lies on classical arguments (as those in \cite[Chap.~VIII, Sec.~1, Prop.~2]{S93}) combined with Cauchy-Schwarz inequality.

\begin{2-THM2} \label{2-THM2}
	Let $p_1$, $p_2$, $\tilde{p}_1$ and $\tilde{p}_2$ be four finite real numbers such that $[p_1, p_2] \subset (\tilde{p}_1, \tilde{p}_2)$. Suppose that $\psi \in \mathcal{C}^2(\R, \R): \R \longrightarrow \R$ is a concave function such that $| \psi' | > 0$ on $[p_1, p_2]$. And assume that $U \in \mathcal{C}^1(\R, \C)$ is a function satisfying
	\begin{equation*}
		supp \, U \subseteq [p_1,p_2] \; .
	\end{equation*}
	Then we have for all $\omega > 0$,
	\begin{align*}
		\left| \int_{\R} U(p) \, e^{i \omega \psi(p)} \, dp \right| \leqslant \Big( C_3(\psi, p_1, p_2) \, \big\| U' \big\|_{L^2(\R)} + C_4(\psi, p_1, p_2) \, \big\| U \big\|_{L^{\infty}(\R)} \Big) \, \omega^{-1} \; .
	\end{align*}
	where
	\begin{align*}
		& \bullet \quad C_3(\psi, p_1, p_2) := (p_2 - p_1)^{\frac{1}{2}} \, \min \Big\{ \big| \psi'(p_1) \big|, \big| \psi'(p_2) \big| \Big\}^{-1} \; ; \\
		& \bullet \quad C_4(\psi, p_1, p_2) := \min \Big\{ \big| \psi'(p_1) \big|, \big| \psi'(p_2) \big| \Big\}^{-1} \; .
	\end{align*}
\end{2-THM2}

\begin{proof}
	Let $\omega > 0$. Since $\psi'$ is monotonic and has a constant sign on $[p_1, p_2]$, we have
	\begin{equation*}
		\forall \, p \in [p_1, p_2] \qquad \big| \psi'(p) \big| \geqslant \min \Big\{ \big| \psi'(p_1) \big|, \big| \psi'(p_2) \big| \Big\} =: m_{p_1, p_2}(\psi') > 0 \; .
	\end{equation*}
	Hence we are allowed to integrate by parts as follows:
	\begin{equation*}
		\int_{\R} U(p) \, e^{i \omega \psi(p)} \, dp = \int_{p_1}^{p_2} U(p) \, e^{i \omega \psi(p)} \, dp = - i \int_{p_1}^{p_2} \left( \frac{U}{\psi'} \right)'\hspace{-1mm}(p) \, e^{i \omega \psi(p)} \, dp \, \omega^{-1} \; .
	\end{equation*}
	Moreover we have
	\begin{align*}
		& \left| - i \int_{p_1}^{p_2} \left( \frac{U}{\psi'} \right)'\hspace{-1mm}(p) \, e^{i \omega \psi(p)} \, dp \right| \\
		& \hspace{1.5cm} \leqslant \left| \int_{p_1}^{p_2} U'(p) \, \psi'(p)^{-1} \, e^{i \omega \psi(p)} \, dp \right| + \int_{p_1}^{p_2} \Big| U(p) \, \psi''(p) \, \psi'(p)^{-2} \Big| \, dp \\
		& \hspace{1.5cm} \leqslant \big\| U' \big\|_{L^2(\R)} \, \big\| (\psi')^{-1} \big\|_{L^2(p_1, p_2)} + \big\| U \big\|_{L^{\infty}(\R)} \int_{p_1}^{p_2} \Big| \psi''(p) \, \psi'(p)^{-2} \Big| \, dp \; ;
	\end{align*}
	as in the preceding proof, we have applied Cauchy-Schwarz inequality to the first integral. Now the hypotheses $\psi'' \leqslant 0$ and $\psi'$ is monotonic with a constant sign allow to carry the following computations out:
	\begin{equation*}
		\int_{p_1}^{p_2} \Big| \psi''(p) \, \psi'(p)^{-2} \Big| \, dp = \left| - \int_{p_1}^{p_2} \psi''(p) \, \psi'(p)^{-2} \, dp \right| = \Big|  \psi'(p_2)^{-1} - \psi'(p_1)^{-1} \Big| \leqslant m_{p_1, p_2}(\psi')^{-1} \; .
	\end{equation*}
	Furthermore, we have
	\begin{equation*}
		\big\| (\psi')^{-1} \big\|_{L^2(p_1, p_2)} \leqslant m_{p_1, p_2}(\psi')^{-1} \, (p_2 - p_1)^{\frac{1}{2}} \; .
	\end{equation*}
	Consequently we obtain
	\begin{align*}
		\left| \int_{\R} U(p) \, e^{i \omega \psi(p)} \, dp \right|
			& \leqslant \Big( m_{p_1, p_2}(\psi')^{-1} \, (p_2 - p_1)^{\frac{1}{2}} \, \big\| U' \big\|_{L^2(\R)} + m_{p_1, p_2}(\psi')^{-1} \, \big\| U \big\|_{L^{\infty}(\R)} \Big) \, \omega^{-1} \; .
	\end{align*}
\end{proof}

\section{Minimization of error estimates and origin of the propagation cone for a family of dispersive equations} \label{sec:applications_de}

We start this section by introducing the Fourier transform $\tf u : \R \longrightarrow \C$ of a function $u : \R \longrightarrow \C$ belonging to the Schwartz space $\mathcal{S}(\R)$:
\begin{equation*}
	\forall \, p \in \R \qquad \tf u(p) := \int_{\R} u(x) \, e^{-ixp} \, dx \; .
\end{equation*}
The Fourier transform defines an invertible operator from $\mathcal{S}(\R)$ onto itself, and can be extended to the space of square-integrable functions $L^2(\R)$ and to the tempered distributions $\mathcal{S}'(\R)$. Moreover, for $u \in L^2(\R)$, Plancherel theorem assures the following equality:
\begin{equation*}
	\forall \, u \in L^2(\R) \qquad \| u \|_{L^2(\R)} = \frac{1}{\sqrt{2\pi}}\big\| \tf u \big \|_{L^2(\R)} \; ;
\end{equation*}
see \cite[Theorem 7.1.6]{H83}.

Consider now a $\mathcal{C}^{\infty}$-function $f : \R \longrightarrow \R$ such that all its derivatives grow at most as a polynomial at infinity and consider the associated operator $f(D) : \mathcal{S}(\R) \longrightarrow \mathcal{S}(\R)$ defined by
\begin{equation*}
	\forall \, x \in \R \qquad f(D) u(x) := \frac{1}{2 \pi} \int_{\R} f(p) \, \tf u(p) \, e^{i x p} \, dp = \tf^{-1} \Big( f \, \tf u \Big)(x) \; ,
\end{equation*}
which can be extended to the tempered distributions $\mathcal{S}'(\R)$. The operator $f(D) : \mathcal{S}'(\R) \longrightarrow \mathcal{S}'(\R)$ is called a \textit{Fourier multiplier} associated to the \emph{symbol} $f$.

Given such an operator, we introduce the following evolution equation on the line,
\begin{equation} \label{eq:evoleq}
	\left\{ \begin{array}{l}
			\left[ i \, \partial_t - f \big(D\big) \right] u_f(t) = 0 \\ [2mm]
			u_f(0) = u_0
	\end{array} \right. \; ,
\end{equation}
for $t \in \R$. If we suppose $u_0 \in \mathcal{S}'(\R)$ then the equation \eqref{eq:evoleq} has a unique solution in $\displaystyle \mathcal{C}^1\big( \R , \mathcal{S}'(\R) \big)$ given by the following solution formula,
\begin{equation*}
	u_f(t) = \tf^{-1} \Big( e^{-i t f} \tf u_0 \Big) \; .
\end{equation*}
We refer to \cite{BAT1994} for a detailed study of this family of equations.

In this this paper, we suppose that the symbol $f$ is strictly convex; an important example of such an equation is given by the free Schrödinger equation whose symbol is $f_S(p) = \frac{1}{2} \, p^2$.\\ For the sake of better presentation of the results, we consider initial data $u_0$ belonging only to the Schwartz space $\mathcal{S}(\R)$ to focus on the approach we propose. We mention that it is possible to extend our results to the case of initial data in $L^2(\R)$ with additional assumptions on regularity and decay; but this falls out of the scope of the paper.\\
Further the initial data are assumed to be in bounded frequency bands, meaning that their Fourier transforms are supported on bounded intervals $[p_1,p_2]$, where $p_1 < p_2$ are two finite real numbers. Under such hypotheses, the solution formula for the equation \eqref{eq:evoleq} defines a function $u_f : \R \times \R \longrightarrow \C$ given by
\begin{equation} \label{eq:formula}
	u_f(t,x) = \frac{1}{2 \pi} \int_{p_1}^{p_2} \tf u_0(p) \, e^{-itf(p) + ixp} \, dp \; .
\end{equation}

We define now the space-time cone related to the symbol $f$ and to the frequency band $[\tilde{p}_1,\tilde{p}_2]$ with origin $(t_0,x_0) \in \R^2$:
\begin{3-DEF2} \label{3-DEF2}
	Let $t_0, x_0, \tilde{p}_1, \tilde{p}_2$ be four finite real numbers such that $\tilde{p}_1 < \tilde{p}_2$ and let $f: \R \longrightarrow \R$ be a symbol.
	\begin{enumerate}
		\item We define the space-time cone $\mathfrak{C}_f\big( [\tilde{p}_1, \tilde{p}_2], (t_0, x_0) \big)$ as follows:
		\begin{equation} \label{eq:cone}
			\mathfrak{C}_f\big( [\tilde{p}_1, \tilde{p}_2], (t_0, x_0) \big) := \left\{ (t,x) \in \big( \R \backslash \{t_0 \} \big) \times \R \, \bigg| \, f'(\tilde{p}_1) \leqslant \frac{x - x_0}{t - t_0} \leqslant f'(\tilde{p}_2) \right\} \; .
		\end{equation}
		\item Let $\mathfrak{C}_f\big( [\tilde{p}_1, \tilde{p}_2], (t_0, x_0) \big)^c$ be the complement of the space-time cone $\mathfrak{C}_f\big( [\tilde{p}_1, \tilde{p}_2], (t_0, x_0) \big)$ in $\big( \R \backslash \{t_0 \} \big) \times \R$ .
	\end{enumerate}	 
\end{3-DEF2}

In this section, we aim at computing time-asymptotic expansions to one term of the solution formula \eqref{eq:formula} for initial data in frequency bands. We show that the resulting first term of these expansions is supported in a space-time cone of type \eqref{eq:cone}, providing asymptotic propagation features for the solutions. In a first step, the origin of the cone is arbitrarily chosen and the remainder estimates are explicit with respect to this origin.
In a second step, we determine the origin of the cone minimizing this remainder estimate.\\

In the following theorem, we provide a time-asymptotic expansion to one term with explicit error estimate of the solution \eqref{eq:formula} in the space-time cone $\mathfrak{C}_f\big( [\tilde{p}_1, \tilde{p}_2], (t_0, x_0) \big)$, where $[p_1, p_2] \subset (\tilde{p}_1, \tilde{p}_2)$ and $(t_0,x_0) \in \R^2$ is arbitrarily chosen. A uniform estimate of \eqref{eq:formula} outside the cone is also established.\\
The proof of Theorem \ref{3-THM1} follows the lines of the one of \cite[Theorem 5.2]{AMD17}: it consists mainly in rewriting the solution formula \eqref{eq:formula} as an oscillatory integral with respect to time and in applying then Theorems \ref{2-THM1} and \ref{2-THM2}. Here the expansion in a cone with arbitrary origin is obtained thanks to a space-time shift in the integral defining \eqref{eq:formula}. And the explicitness of the remainder with respect to the origin is possible thanks to the new remainder estimate given in Theorem \ref{2-THM1}, allowing the application of Plancherel theorem.

\begin{3-THM1} \label{3-THM1}
	Let $p_1$, $p_2$, $\tilde{p}_1$ and $\tilde{p}_2$ be four finite real numbers such that $[p_1, p_2] \subset (\tilde{p}_1, \tilde{p}_2)$. Suppose that $u_0 \in \mathcal{S}(\R)$ is a function whose Fourier transform satisfies
	\begin{equation*}
		supp \, \tf u_0 \subseteq [p_1,p_2] \; .
	\end{equation*}
	Fix $(t_0,x_0) \in \R^2$. Then
	\begin{enumerate}
		\item for all $(t,x) \in \mathfrak{C}_f\big( [\tilde{p}_1, \tilde{p}_2], (t_0, x_0) \big)$, we have
		\begin{align}
			& \left| u_f(t, x) - \frac{1}{\sqrt{2 \pi}} \, e^{- sgn(t-t_0) i \frac{\pi}{4}} \, e^{-itf(p_0(t,x)) + ix p_0(t,x)} \, \frac{\tf u_0 \big( p_0(t,x) \big)}{\sqrt{f''\big( p_0(t,x) \big)}} \, |t-t_0|^{-\frac{1}{2}} \right| \nonumber \\
			 & \hspace{2cm} \leqslant \Big( C_5(f, \delta, \tilde{p}_1, \tilde{p}_2) \, \big\| (.-x_0) \, u_f(t_0,.) \big\|_{L^2(\R)} \nonumber \\
			 & \hspace{4cm} + C_6(f, \delta, \tilde{p}_1, \tilde{p}_2) \, \big\| u_0 \big\|_{L^1(\R)} \Big) \, |t-t_0|^{- \delta} \; , \label{eq:remainder_est1}
		\end{align}
		where the real number $\delta$ is arbitrarily chosen in $\big( \frac{1}{2}, \frac{3}{4} \big)$ and
		\begin{align*}
			& \bullet \quad p_0(t,x) := (f')^{-1} \hspace{-1mm} \left( \frac{x-x_0}{t-t_0} \right) \; ; \\
			& \bullet \quad C_5(\psi, \delta, \tilde{p}_1, \tilde{p}_2) := \frac{2^{\delta+\frac{1}{2}} \, L(\delta)}{\sqrt{\pi} \sqrt{3-4\delta}} \, \big( \tilde{p}_2 - \tilde{p}_1 \big)^{\frac{3 - 4 \delta}{2}} \, c_5(f, \delta, \tilde{p}_1, \tilde{p}_2) \; ; \\
			& \bullet \quad C_6(\psi, \delta, \tilde{p}_1, \tilde{p}_2) := \frac{2^{\delta - 2} L(\delta)}{\pi(1 - \delta)} \, \big( \tilde{p}_2 - \tilde{p}_1 \big)^{2-2\delta} \, c_6(f, \delta, \tilde{p}_1, \tilde{p}_2) \; ; \\
			& \bullet \quad c_5(f, \delta, \tilde{p}_1, \tilde{p}_2) := \big\| f'' \big\|_{L^{\infty}(\tilde{p}_1, \tilde{p}_2)}^{\frac{3}{2}-\delta} \, \min_{[\tilde{p}_1, \tilde{p}_2]}\big\{f''\big\}^{-\frac{3}{2}} \; ; \\
			&  \bullet \quad c_6(\psi, \delta, \tilde{p}_1, \tilde{p}_2) := \big\| f'' \big\|_{L^{\infty}(\tilde{p}_1, \tilde{p}_2)}^{\frac{5}{2}-\delta} \, \big\| f^{(3)} \big\|_{L^{\infty}(\tilde{p}_1, \tilde{p}_2)} \, \min_{[\tilde{p}_1, \tilde{p}_2]} \big\{f''\big\}^{-\frac{7}{2}} \\
				& \hspace{4.5cm} + \frac{1}{3} \, \big\| f'' \big\|_{L^{\infty}(\tilde{p}_1, \tilde{p}_2)}^{\frac{7}{2}-\delta} \, \big\| f^{(3)} \big\|_{L^{\infty}(\tilde{p}_1, \tilde{p}_2)} \, \min_{[\tilde{p}_1, \tilde{p}_2]} \big\{f''\big\}^{-\frac{9}{2}} \; .
		\end{align*}		
		The constant $L(\delta) > 0$ is defined in Lemma \ref{2-LEM1} iii);
		\item for all $(t,x) \in \mathfrak{C}_f\big( [\tilde{p}_1, \tilde{p}_2], (t_0, x_0) \big)^c$, we have
		\begin{align}
			\big| u_f(t,x) \big|
				& \leqslant \Big( C_7(f, p_1, p_2, \tilde{p}_1, \tilde{p}_2) \, \big\| (.-x_0) \, u_f(t_0,.) \big\|_{L^2(\R)} \nonumber \\
				& \hspace{1.5cm} +  C_8(f, p_1, p_2, \tilde{p}_1, \tilde{p}_2) \, \big\| u_0 \big\|_{L^1(\R)} \Big) \, |t-t_0|^{-1} \; , \label{eq:remainder_est2}
		\end{align}
		where
		\begin{align*}
			& \bullet \quad C_7(f, p_1, p_2, \tilde{p}_1, \tilde{p}_2) := \frac{1}{\sqrt{2\pi}} \, (p_2 - p_1)^{\frac{1}{2}}\, \min \hspace{-1mm} \big\{ f'(p_1) - f'(\tilde{p}_1), f'(\tilde{p}_2) - f'(p_2) \big\}^{-1} \; ; \\
			& \bullet \quad C_8(f, p_1, p_2, \tilde{p}_1, \tilde{p}_2) := \frac{1}{2\pi} \, \min \hspace{-1mm} \big\{ f'(p_1) - f'(\tilde{p}_1), f'(\tilde{p}_2) - f'(p_2) \big\}^{-1} \; .
		\end{align*}
	\end{enumerate}
\end{3-THM1}

\begin{proof}
	The cases where $t > t_0$ and $t < t_0$ are distinguished for the sake of readability.\\
	
	\underline{Case 1:} $t > t_0$\\
	We rewrite the solution formula as an oscillatory integral by proceeding as follows\footnote{In \cite{AMD17, D17-1}, the parameters $t_0$ and $x_0$ are implicitly equal to $0$. Allowing these parameters to be arbitrary produces a space-time shift in the solution formula and permits to consider cones with arbitrary origin.}
	\begin{align*}
		u_f(t,x)	& = \frac{1}{2 \pi} \int_{\R} \tf u_0(p) \, e^{-i t f(p) + i x p} \, dp \\
					& = \int_{\R} \frac{1}{2 \pi} \, \tf u_0(p) \, e^{-i t_0 f(p) + i x_0 p} \, e^{i (t-t_0) \big( \frac{x-x_0}{t-t_0} p \, - \, f(p) \big)} \, dp \\
					& = \int_{\R} \mathbf{U}_f(p,t_0,x_0) \, e^{i (t-t_0) \Psi_f(p,t,x,t_0,x_0)} \, dp \\
					& =: I_f(t, x, u_0, t_0, x_0) \; .
	\end{align*}
	We note that the amplitude
	\begin{equation*}
		\mathbf{U}_f(p,t_0,x_0) := \frac{1}{2 \pi} \, \tf u_0(p) \, e^{-i t_0 f(p) + i x_0 p} \; ,
	\end{equation*}
	which is actually the Fourier transform of $\frac{1}{2 \pi} u_f(t_0, \, . + x_0)$, is a $\mathcal{C}^\infty$-function (with respect to the variable $p$) whose support is included in $[p_1,p_2]$. The phase function
	\begin{equation*}
		\Psi_f(p,t,x,t_0,x_0) := \frac{x-x_0}{t-t_0} p \, - \, f(p)
	\end{equation*}
	is a $\mathcal{C}^{\infty}$-function on $\R$ which is strictly concave since we have supposed $f'' > 0$ in this section.\\
	Now we remark that the existence of a stationary point for the phase inside the interval $\tilde{I} := (\tilde{p}_1, \tilde{p}_2)$ depends on the value of $\frac{x-x_0}{t-t_0}$: it exists and is unique if and only if $\frac{x-x_0}{t-t_0} \in f'\big( \tilde{I}\big)$. In this case, the stationary point $p_0(t,x)$ is given by
	\begin{equation*}
		p_0(t,x) = (f')^{-1} \hspace{-1mm} \left( \frac{x-x_0}{t-t_0} \right) \; .
	\end{equation*}
	Let us now distinguish two sub-cases to apply Theorem \ref{2-THM1} and Theorem \ref{2-THM2}.
	\begin{itemize}
		\item \emph{Case} $\frac{x-x_0}{t-t_0} \in f'\big( \tilde{I}\big)$. In this case, the stationary point belongs to $\tilde{I}$. Hence we are allowed to apply Theorem \ref{2-THM1} to the oscillatory integral $I_f(t, x, u_0, t_0, x_0)$ with $\omega = t - t_0$:
		\begin{align*}
			& \left| I_f(t, x, u_0, t_0, x_0) - \frac{1}{\sqrt{2 \pi}} \, e^{-i \frac{\pi}{4}} \, e^{-itf(p_0(t,x)) + ix p_0(t,x)} \, \frac{\tf u_0 \big( p_0(t,x) \big)}{\sqrt{f''\big( p_0(t,x) \big)}} \, (t-t_0)^{-\frac{1}{2}} \right| \\
			 & \hspace{2cm} \leqslant \frac{1}{2 \pi} \, \bigg( C_1(\Psi_f, \delta, \tilde{p}_1, \tilde{p}_2) \, \Big\| \partial_p \left[ \tf u_0(\cdot) \, e^{- i t_0 f(\cdot) + i x_0 \cdot} \right] \Big\|_{L^2(\R)} \\
			 & \hspace{4cm} + C_2(\Psi_f, \delta, \tilde{p}_1, \tilde{p}_2) \, \big\| \tf u_0 \big\|_{L^{\infty}(\R)} \bigg) \, (t-t_0)^{- \delta} \; ,
		\end{align*}
		with $\delta \in \big( \frac{1}{2}, \frac{3}{4} \big)$ and the constants $C_1(\Psi_f, \delta, \tilde{p}_1, \tilde{p}_2)$, $C_2(\Psi_f, \delta, \tilde{p}_1, \tilde{p}_2) > 0$ are defined in Theorem \ref{2-THM1}. Since we have
		\begin{equation*}
			\partial_p^{(2)} \Psi_f(p,t,x,t_0,x_0) = -f''(p) \quad , \quad \partial_p^{(3)} \Psi_f(p,t,x,t_0,x_0) = -f^{(3)}(p) \; ,
		\end{equation*}
		and since the constants $C_1(\Psi_f, \delta, \tilde{p}_1, \tilde{p}_2)$ and $C_2(\Psi_f, \delta, \tilde{p}_1, \tilde{p}_2)$ depend only on the se\-cond and third derivatives (with respect to $p$) of the phase, we can claim that these constants depend on $f$ rather than $\Psi_f$. Furthermore, Plancherel theorem and standard properties of the Fourier transform provide
		\begin{align*}
			\Big\| \partial_p \left[ \tf u_0(.) e^{-it f(\cdot) + i x_0 \cdot} \right] \Big\|_{L^2(\R)}	& = \sqrt{2\pi} \, \Big\| x \longmapsto x \, \tf^{-1}\big[\tf u_0(.) e^{-it_0 f(\cdot) + i x_0 \cdot} \big](x) \Big\|_{L^2(\R)} \\
							& = \sqrt{2\pi} \, \Big\| x \longmapsto x \, \tf^{-1}\big[\tf u_0(.) \, e^{-it_0 f(\cdot)} \big](x + x_0) \Big\|_{L^2(\R)} \\
							& = \sqrt{2\pi} \, \Big\| x \longmapsto (x-x_0) \, \tf^{-1}\big[\tf u_0(.) \, e^{-it_0 f(\cdot)} \big](x) \Big\|_{L^2(\R)} \\
							& = \sqrt{2\pi} \, \Big\| x \longmapsto (x-x_0) \, u_f(t_0,x) \Big\|_{L^2(\R)} \; .
		\end{align*}
		Hence we obtain finally
		\begin{align*}
			& \left| u_f(t, x) - \frac{1}{\sqrt{2 \pi}} \, e^{-i \frac{\pi}{4}} \, e^{-itf(p_0(t,x)) + ix p_0(t,x)} \, \frac{\tf u_0 \big( p_0(t,x) \big)}{\sqrt{f''\big( p_0(t,x) \big)}} \, (t-t_0)^{-\frac{1}{2}} \right| \\
			 & \hspace{2cm} \leqslant \Bigg( \frac{1}{\sqrt{2 \pi}} \, C_1(-f, \delta, \tilde{p}_1, \tilde{p}_2) \, \big\| (.-x_0) \, u_f(t_0,.) \big\|_{L^2(\R)} \\
			 & \hspace{4cm} + \frac{1}{2 \pi} \, C_2(-f, \delta, \tilde{p}_1, \tilde{p}_2) \, \big\| u_0 \big\|_{L^1(\R)} \Bigg) \, (t-t_0)^{- \delta} \; ,
		\end{align*}		
		where we have used the classical estimate $\big\| \tf u_0 \big\|_{L^\infty(\R)} \leqslant \| u_0 \|_{L^1(\R)}$.
		\item \emph{Case} $\frac{x-x_0}{t-t_0} \notin f'\big(\tilde{I}\big)$. As previously, we rewrite the solution formula as the oscillatory integral $I_f(t,x,u_0,t_0,x_0)$. Here the phase $\Psi_f(.,t,x,t_0,x_0)$ has no stationary point inside the interval $\tilde{I} = (\tilde{p}_1, \tilde{p}_2)$ and one has
		\begin{equation*}
			\forall \, p \in [p_1, p_2] \qquad \Big| \partial_p \Psi_f(p,t,x,t_0,x_0) \Big| = \left| \frac{x-x_0}{t-t_0} - f'(p) \right| \geqslant m_{\tilde{I}}(f) > 0
		\end{equation*}
		where $m_{\tilde{I}}(f) := \min \big\{ f'(p_1) - f'(\tilde{p}_1), f'(\tilde{p}_2) - f'(p_2) \big\}$. Consequently we can apply Theorem \ref{2-THM2} which provides
		\begin{align*}
			\big| u_f(t,x) \big|
				& \leqslant \Bigg( \frac{1}{\sqrt{2 \pi}} \, C_3(-f, \tilde{p}_1, \tilde{p}_2) \, \big\| (.-x_0) \, u_f(t_0,.) \big\|_{L^2(\R)} \\
				& \hspace{1.5cm} + \frac{1}{2 \pi} \, C_4(-f, \tilde{p}_1, \tilde{p}_2) \, \big\| u_0 \big\|_{L^1(\R)} \Bigg) (t-t_0)^{-1} \; ,
		\end{align*}
		where the constants $C_3(-f, \tilde{p}_1, \tilde{p}_2)$, $C_4(-f, \tilde{p}_1, \tilde{p}_2) > 0$ are defined in Theorem \ref{2-THM2}.
	\end{itemize}	 
	
	\underline{Case 2:} $t < t_0$\\
	Here we have
	\begin{align*}
		u_f(t,x)	& = \int_{\R} \frac{1}{2 \pi} \, \tf u_0(p) \, e^{-i t_0 f(p) + i x_0 p} \, e^{i (t-t_0) \big( \frac{x-x_0}{t-t_0} p \, - \, f(p) \big)} \, dp \\
					& = \overline{\int_{\R} \frac{1}{2 \pi} \, \overline{\tf u_0(p) \, e^{-i t_0 f(p) + i x_0 p}} \, e^{i (t_0-t) \big( \frac{x-x_0}{t-t_0} p \, - \, f(p) \big)} \, dp} \\
					& = \overline{\int_{\R} \overline{\mathbf{U}_f(p,t_0,x_0)} \, e^{i (t_0-t) \Psi_f(p,t,x,t_0,x_0)} \, dp} \; .
	\end{align*}
	Following the arguments and computations of the preceding case $t > t_0$, we obtain
	\begin{align*}
		& \left| u_f(t, x) - \frac{1}{\sqrt{2 \pi}} \, e^{+i \frac{\pi}{4}} \, e^{-itf(p_0(t,x)) + ix p_0(t,x)} \, \frac{\tf u_0 \big( p_0(t,x) \big)}{\sqrt{f''\big( p_0(t,x) \big)}} \, (t_0-t)^{-\frac{1}{2}} \right| \\
	    & \hspace{2cm} \leqslant \bigg( \frac{1}{\sqrt{2 \pi}} \, C_1(-f, \delta, \tilde{p}_1, \tilde{p}_2) \big\| (.-x_0) \, u_f(t_0,.) \big\|_{L^2(\R)} \\
		& \hspace{4cm} + \frac{1}{2\pi} \, C_2(-f, \delta, \tilde{p}_1, \tilde{p}_2) \big\| u_0 \big\|_{L^1(\R)} \bigg) \, (t_0-t)^{- \delta} \; ,
	\end{align*}
	for all $(t,x) \in (-\infty, t_0) \times \R$ such that $\frac{x-x_0}{t-t_0} \in f' \big( \tilde{I} \big)$, and
	\begin{align*}
		\big| u_f(t,x) \big|
			& \leqslant \Bigg( \frac{1}{\sqrt{2 \pi}} \, C_3(-f, \tilde{p}_1, \tilde{p}_2) \, \big\| (.-x_0) \, u_f(t_0,.) \big\|_{L^2(\R)} \\
			& \hspace{1.5cm} + \frac{1}{2 \pi} \, C_4(-f, \tilde{p}_1, \tilde{p}_2) \, \big\| u_0 \big\|_{L^1(\R)} \Bigg) (t_0-t)^{-1} \; ,
	\end{align*}
	for all $(t,x) \in (-\infty, t_0) \times \R$ such that $\frac{x-x_0}{t-t_0} \notin f' \big( \tilde{I} \big)$, leading to the desired estimates.
\end{proof}

We define now the following moment-type and variance-type quantities for normalized $u \in L^2(\R)$.

\begin{3-DEF1} \label{3-DEF1}
	Choose $u \in L^2(\R)$ such that $\| u \|_{L^2(\R)} = 1$ and let $f : \R \longrightarrow \R$ be a symbol. If they exist, we define the real numbers $\mathcal{M}_f(u)$ and $\mathcal{V}_f(u)$ as follows,
	\begin{align*}
		\mathcal{M}_f(u) := \int_{\R} f(x) \, \big| u(x) \big|^2 \, dx \quad , \quad \mathcal{V}_f(u) := \mathcal{M}_{f^2}(u) - \mathcal{M}_f(u)^2 \; .
	\end{align*}
	If $f(x) = x$, then we note for simplicity
	\begin{equation*}
		\mathcal{M}_1(u) := \mathcal{M}_f(u) \quad , \quad \mathcal{M}_2(u) := \mathcal{M}_{f^2}(u) \quad , \quad \mathcal{V}(u) := \mathcal{V}_f(u) \; .
	\end{equation*}
\end{3-DEF1}

\begin{3-REM0} \em
	The above quantities $\mathcal{M}_1(u)$, $\mathcal{M}_2(u)$ and $\mathcal{V}(u)$ are respectively the mean, the second moment and the variance of $|u|^2$.
\end{3-REM0}

The following lemma will permit to determine the space-time cone in which the remainder bound given in Theorem \ref{3-THM1} is minimal with respect to $(t_0,x_0)$; see Corollary \ref{3-COR1}. This is based on the minimisation of the function $(t_0, x_0) \longmapsto \big\| (.-x_0) u_f(t_0,.) \big\|_{L^2(\R)}^2$ which is the moment of order 2 of $x \longmapsto \big| u_f(t_0,x+x_0) \big|^2$ for fixed $(t_0,x_0) \in \R^2$.

\begin{3-LEM1} \label{3-LEM1}
	Suppose that the hypotheses of Theorem \ref{3-THM1} are satisfied and suppose in addition that $\| u_0 \|_{L^2(\R)} = 1$.	Then the function $g : \R^2 \longrightarrow \R_+$ defined by
	\begin{equation*}
		g(t_0,x_0) = \big\| (.-x_0) \, u_f(t_0,.) \big\|_{L^2(\R)}^2
	\end{equation*}
	has a global minimum at $(t^*, x^*) \in \R^2$ with
	\begin{align*}
		& \bullet \quad t^* = \argmin_{\tau \in \R} \mathcal{V}\big(u_f(\tau,.)\big) \; ; \\
		& \bullet \quad x^* = \mathcal{M}_1\big( u_f(t^*,.) \big) \; .
	\end{align*}
\end{3-LEM1}

\begin{proof}
	For fixed $t_0 \in \R$, differentiating twice the function $g(t_0, .)$ with respect to its second argument shows that
	\begin{equation*}
		\partial_{x_0}^{(2)} g(t_0, x_0) = 2 > 0 \; ,
	\end{equation*}
	Hence $g(t_0,.)$ is a polynomial function of degree 2 whose unique global minimum is
	\begin{equation*}
		\tilde{x}(t_0) = \int_\R x \big| u_f(t_0,x) \big|^2 \, dx = \mathcal{M}_1\big(u_f(t_0,.) \big) \; .
	\end{equation*}
	It follows that
	\begin{equation*}
		g\big(t_0, \tilde{x}(t_0) \big) = \int_\R \Big(x-\mathcal{M}_1\big(u_f(t_0,.) \big) \Big)^2 \, \big| u_f(t_0,x) \big|^2 \, dx = \mathcal{V}\big(u_f(t_0,.) \big) \; .
	\end{equation*}
	Lemma \ref{A-PROP2} assures that $ t_0 \in \R \longmapsto \mathcal{V}\big(u_f(t_0,.) \big) \in \R_+$ is a polynomial of degree 2 whose leading coefficient is $\mathcal{V}_{f'}\big( \frac{1}{\sqrt{2 \pi}} \tf u_0 \big)$. Since we have by simple calculations,
	\begin{align*}
		\mathcal{V}_{f'} \hspace{-1mm} \left( \frac{1}{\sqrt{2 \pi}} \tf u_0 \right)
			& = \mathcal{M}_{f'^{\, 2}} \hspace{-1mm} \left( \frac{1}{\sqrt{2 \pi}} \, \tf u_0 \right) - \mathcal{M}_{f'} \hspace{-1mm} \left( \frac{1}{\sqrt{2 \pi}} \, \tf u_0 \right)^2 \\
			& = \frac{1}{2 \pi} \int_\R \left( f'(p) - \mathcal{M}_{f'} \hspace{-1mm} \left( \frac{1}{\sqrt{2 \pi}} \, \tf u_0 \right) \right)^2 \big| \tf u_0(p) \big|^2 \, dp \; ,
	\end{align*}
	and since $f$ is supposed to be strictly convex in this paper, the leading coefficient $\mathcal{V}_{f'}\big( \frac{1}{\sqrt{2 \pi}} \tf u_0 \big)$ is necessarily positive. Thus the function $t_0 \in \R \longmapsto g\big(t_0, \tilde{x}(t_0) \big) \in \R_+$ has a global minimum at a certain $t^* \in \R$, \emph{i.e.},
	\begin{equation*}
		t^* = \argmin_{\tau \in \R} \mathcal{V}\big(u_f(\tau,.)\big) \; .
	\end{equation*}
	Finally we define
	\begin{equation*}
		x^* := \tilde{x}(t^*) = \mathcal{M}_1\big(u_f(t^*,.) \big) \; .
	\end{equation*}
\end{proof}

\begin{3-REM1} \label{3-REM1} \em
	The polynomial nature of the function $t_0 \in \R \longmapsto g\big(t_0, \tilde{x}(t_0) \big) \in \R_+$ permits to derive the following formula for $t^*$:
	\begin{align}
		t^*
			& = \frac{1}{\mathcal{V}_{f'}\big( \frac{1}{\sqrt{2 \pi}} \, \tf u_0 \big)} \left( - \frac{1}{2 \pi} \Im \bigg(  \int_{\R} f'(p) \, \tf u_0(p) \, \overline{\big( \tf u_0 \big)'(p)} \, dp \right) \nonumber \\
			& \hspace{1.5cm} + \mathcal{M}_{f'} \hspace{-1mm} \left( \frac{1}{\sqrt{2 \pi} } \, \tf u_0 \right) \mathcal{M}_1(u_0) \bigg) \; . \label{eq:formula_t}
	\end{align}
	Furthermore, from Lemma \ref{A-PROP1}, we have
	\begin{equation*}
		\mathcal{M}_1\big(u_f(t^*,.) \big) = \mathcal{M}_{f'} \hspace{-1mm} \left( \frac{1}{\sqrt{2 \pi}} \, \tf u_0 \right) t^* + \mathcal{M}_1(u_0) \; ;
	\end{equation*}
	inserting formula \eqref{eq:formula_t} into the preceding equality provides
	\begin{align*}
		x^* 
			& = \frac{1}{\mathcal{V}_{f'}\big( \frac{1}{\sqrt{2 \pi}} \, \tf u_0 \big)} \bigg( - \frac{1}{2 \pi} \Im \hspace{-1mm} \left(  \int_{\R} f'(p) \, \tf u_0(p) \, \overline{\big( \tf u_0 \big)'(p)} \, dp \right) \mathcal{M}_{f'} \hspace{-1mm} \left( \frac{1}{\sqrt{2 \pi} } \, \tf u_0 \right) \\
			& \hspace{1.5cm} + \mathcal{M}_{f'^2} \hspace{-1mm} \left( \frac{1}{\sqrt{2 \pi} } \, \tf u_0 \right) \mathcal{M}_1(u_0) \bigg) \; .
	\end{align*}
\end{3-REM1}

The final result of this section is a direct consequence of Theorem \ref{3-THM1} and of Lemma \ref{3-LEM1}: it shows that the bounds of the error estimates appearing in Theorem \ref{3-THM1} are minimised by putting the origin of the cone at the point $\big( t^*, x^* \big)$ defined above.

\begin{3-COR1} \label{3-COR1}
	Suppose that the hypotheses of Theorem \ref{3-THM1} are satisfied and suppose in addition that $\| u_0 \|_{L^2(\R)} = 1$. Then the $(t_0, x_0)$-dependent right-hand sides of estimates \eqref{eq:remainder_est1} and \eqref{eq:remainder_est2} in Theorem \ref{3-THM1} have a global minimum at the point $\big( t^*, x^* \big) \in \R^2$ with
	\begin{align*}
		& \bullet \quad t^* = \argmin_{\tau \in \R} \mathcal{V}\big(u_f(\tau,.)\big) \; ; \\
		& \bullet \quad x^* = \mathcal{M}_1\big( u_f(t^*,.) \big) \; .
	\end{align*}
	In this case, for all $(t,x) \in \mathfrak{C}_f\big( [\tilde{p}_1, \tilde{p}_2], (t^*, x^*) \big)$, we have
	\begin{align*}
		& \left| u_f(t, x) - \frac{1}{\sqrt{2 \pi}} \, e^{- sgn(t-t^*) i \frac{\pi}{4}} \, e^{-itf(p_0(t,x)) + ix p_0(t,x)} \, \frac{\tf u_0 \big( p_0(t,x) \big)}{\sqrt{f''\big( p_0(t,x) \big)}} \, |t-t^*|^{-\frac{1}{2}} \right| \nonumber \\
		 & \hspace{2cm} \leqslant \Big( C_5(f, \delta, \tilde{p}_1, \tilde{p}_2) \, \min_{\tau \in \R} \left(\sqrt{\mathcal{V}\big(u_f(\tau,.) \big)}\right) \nonumber \\
		 & \hspace{4cm} + C_6(f, \delta, \tilde{p}_1, \tilde{p}_2) \, \big\| u_0 \big\|_{L^1(\R)} \Big) \, |t-t^*|^{- \delta} \; ,
	\end{align*}
	where the real number $\delta$ is arbitrarily chosen in $\big( \frac{1}{2}, \frac{3}{4} \big)$ and $p_0(t,x) := (f')^{-1}\big( \frac{x-x^*}{t - t^*} \big)$, and for all $(t,x) \in \mathfrak{C}_f\big( [\tilde{p}_1, \tilde{p}_2], (t^*, x^*) \big)^c$, we have
	\begin{align*}
		\big| u_f(t,x) \big|
			& \leqslant \Big( C_7(f, p_1, p_2, \tilde{p}_1, \tilde{p}_2) \, \min_{\tau \in \R} \left(\sqrt{\mathcal{V}\big(u_f(\tau,.) \big)}\right) \nonumber \\
			& \hspace{1.5cm} +  C_8(f, p_1, p_2, \tilde{p}_1, \tilde{p}_2) \, \big\| u_0 \big\|_{L^1(\R)} \Big) \, |t-t^*|^{-1} \; ; \label{eq:remainder_est2}
	\end{align*}
	all the constants are defined in Theorem \ref{3-THM1}.
\end{3-COR1}

\begin{proof}
	For fixed initial datum and $\delta \in \big( \frac{1}{2}, \frac{3}{4} \big)$, it is clear that minimizing the remainder bounds \eqref{eq:remainder_est1} and \eqref{eq:remainder_est2} with respect to $(t_0, x_0)$ is equivalent to minimizing the function $g : \R^2 \longrightarrow \R_+$ defined in Lemma \ref{3-LEM1}. This lemma affirms that $g$ has a global minimum at $(t^*, x^*)$ defined above. Furthermore we have
	\begin{align*}
		\big\| (.-x^*) \, u_f(t^*,.) \big\|_{L^2(\R)}^2
			& = \int_\R \Big(x-\mathcal{M}_1\big(u_f(t^*,.) \big) \Big)^2 \, \big| u_f(t^*,x) \big|^2 \, dx \\
			& = \mathcal{V}\big(u_f(t^*,.) \big) \; ,
	\end{align*}
	which is equal to the minimum of $\tau \in \R \longmapsto \mathcal{V}\big(u_f(\tau,.) \big)$ by the definition of $t^*$. This ends the proof.
\end{proof}

\section{Mean position and variance stable under time-asymp\-to\-tic approximations} \label{sec:preservation}

Here we are interested in the mean position and the variance of the first term of the expansions given in Theorem \ref{3-THM1}. We exploit the flexibility inherited from the preceding section to choose the origin of the space-time cone, in which we expand the solution of equation \eqref{eq:evoleq}, in such a way that the associated first term and the solution share the same mean position and the difference between the variances is an explicit constant. It turns out that such a cone corresponds to the one in which the $(t_0,x_0)$-dependent remainder estimate from Theorem \ref{3-THM1} is minimized.\\
This section illustrates that the refined method we have developed in the present paper offers approximations of the solution of equation \eqref{eq:evoleq} describing precisely its time-asymptotic propagation features.\\

Let $u_0 \in \mathcal{S}(\R)$ such that
\begin{equation*}
	supp \, \tf u_0 \subseteq [p_1, p_2] \; ,
\end{equation*}
where $p_1 < p_2$ are two finite real numbers, let $f : \R \longrightarrow \R$ be a strictly convex symbol and choose $(t_0, x_0) \in \R^2$. We define $H_f(.,.,u_0,t_0,x_0) : \big( \R \backslash \{t_0 \}\big) \times \R \longrightarrow \C$ as
\begin{equation*} \label{eq:H_f}
	H_f(t,x,u_0,t_0,x_0) := \frac{1}{\sqrt{2 \pi}} \, e^{- sgn(t - t_0) i \frac{\pi}{4}} \, e^{-it f( p_0(t,x)) + ix p_0(t,x)} \, \frac{\tf u_0 \big(p_0(t,x) \big)}{\sqrt{f''\big(p_0(t,x) \big)}} \, |t - t_0|^{-\frac{1}{2}} \; ,
\end{equation*}
with $p_0(t,x) := (f')^{-1}\big( \frac{x-x_0}{t-t_0} \big)$. The function $H_f(.,.,u_0,t_0,x_0)$ is actually the first term of the expansion given in Theorem \ref{3-THM1}; we recall that it is supported in the cone $\mathfrak{C}_f \big( [p_1,p_2], (t_0,x_0) \big)$.\\
In the following proposition, we compute the mean position of $H_f(t,.,u_0,t_0,x_0)$ for all $t \neq t_0$.

\begin{4-PROP1} \label{4-PROP1}
	Let $(t_0,x_0) \in \R^2$. Suppose that the hypotheses of Theorem \ref{3-THM1} are satisfied and suppose in addition that $\| u_0 \|_{L^2(\R)} = 1$. Then for all $t \in \R \backslash \{ t_0 \}$, we have
	\begin{equation*}
		\mathcal{M}_1\Big( H_f(t,.,u_0,t_0,x_0) \Big) = x_0 + \mathcal{M}_{f'} \hspace{-1mm} \left( \frac{1}{\sqrt{2 \pi}} \, \tf u_0 \right) (t-t_0) \; .
	\end{equation*}
\end{4-PROP1}

\begin{proof}
	Let $t \in \R \backslash \{t_0\}$. First of all, we note that
	\begin{equation*}
		x_0 + f'(p_1)(t-t_0) < x_0 + f'(p_2)(t-t_0) \quad \Longleftrightarrow \quad t > t_0 \; .
	\end{equation*}
	Hence using the definition of $H_f(t,x,u_0,t_0,x_0)$ given just above, we have
	\begin{align*}
		\int_{\R} x \, \Big|H_f(t,x,u_0,t_0,x_0) \Big|^2 \, dx
			& = \frac{sgn(t-t_0)}{2 \pi} \int_{x_0 + f'(p_1) (t-t_0)}^{x_0 + f'(p_2) (t-t_0)} x \, \frac{ \big| \tf u_0 \big( p_0(t,x) \big) \big|^2}{f''\big(p_0(t,x) \big)} \, dx \, |t - t_0 |^{-1} \\
			& = \frac{1}{2 \pi} \int_{x_0 + f'(p_1) (t-t_0)}^{x_0 + f'(p_2) (t-t_0)} x \, \frac{ \big| \tf u_0 \big( p_0(t,x) \big) \big|^2}{f''\big(p_0(t,x) \big)} \, dx \, (t - t_0 )^{-1} \; .
	\end{align*}
	We make now the change of variable $x = x_0 + f'(p)(t-t_0)$ to obtain the desired result:
	\begin{align*}
		\int_{\R} x \, \Big| H_f \big(t,x, u_0, t_0, x_0 \big) \Big|^2 \, dx	& = \frac{1}{2 \pi} \int_{p_1}^{p_2} \big(x_0 + f'(p) (t-t_0 ) \big) \big| \tf u_0 (p) \big|^2 \, dp \\
		& = x_0 + \mathcal{M}_{f'} \hspace{-1mm} \left( \frac{1}{\sqrt{2 \pi}} \tf u_0 \right) \, t - \mathcal{M}_{f'} \hspace{-1mm} \left( \frac{1}{\sqrt{2 \pi}} \tf u_0 \right) \, t_0 \; ;
	\end{align*}
	note that we have used the fact that $\frac{1}{2 \pi} \big\| \tf u_0 \big\|_{L^2(p_1,p_2)}^2 = 1$, which is a direct consequence of the assumption $\| u_0 \|_{L^2(\R)} = 1$.
\end{proof}

In the following result, we prove that the mean positions of the solution $u_f(t,.)$ of equation \eqref{eq:evoleq} and of the first term $H_f(t,.,u_0,t_0,x_0)$ are equal if and only if $(t_0,x_0)$ belongs to the space-time line given by the mean position of the solution of \eqref{eq:evoleq}.

\begin{4-PROP2} \label{4-PROP2}
	Let $(t_0,x_0) \in \R^2$. Suppose that the hypotheses of Theorem \ref{3-THM1} are satisfied and suppose in addition that $\| u_0 \|_{L^2(\R)} = 1$. Then for all $t \in \R \backslash \{ t_0 \}$, we have
	\begin{equation*}
		\mathcal{M}_1 \big( u_f (t,.) \big) = \mathcal{M}_1\big( H_f(t,.,u_0,t_0,x_0) \big) \quad \Longleftrightarrow \quad x_0 = \mathcal{M}_1 \big( u_f (t_0,.) \big) \; .
	\end{equation*}
\end{4-PROP2}

\begin{proof}
	Let $t \in \R \backslash \{t_0\}$. According to Propositions \ref{4-PROP1} and \ref{A-PROP1}, we have
	\begin{align*}
		& \bullet \quad \mathcal{M}_1\big( u_f(t,.) \big) = \mathcal{M}_1(u_0) + \mathcal{M}_{f'} \hspace{-1mm} \left( \frac{1}{\sqrt{2 \pi}} \tf u_0 \right) t \; ; \\
		& \bullet \quad \mathcal{M}_1\big( H_f(t,.,u_0,t_0,x_0) \big) = x_0 + \mathcal{M}_{f'} \hspace{-1mm} \left( \frac{1}{\sqrt{2 \pi}} \, \tf u_0 \right) (t-t_0) \; .
	\end{align*}
	Hence these two mean positions are equal if and only if
	\begin{equation*}
		\mathcal{M}_1(u_0) = x_0 - \mathcal{M}_{f'} \hspace{-1mm} \left( \frac{1}{\sqrt{2 \pi}} \, \tf u_0 \right) t_0 \; ,
	\end{equation*}
	which is equivalent to $x_0 = \mathcal{M}_1 \big( u_f (t_0,.) \big)$.
\end{proof}

\begin{4-REM1} \label{4-REM1} \em
	The definition of $x^*$ from Lemma \ref{3-LEM1} (or Corollary \ref{3-COR1}) and the preceding proposition assure that the mean positions of $u_f(t,.)$ and $H_f(t,.,u_0, t^*,x^*)$ are equal for all $t \neq t_0$.
\end{4-REM1}

In the two following results, we focus on the variances of the solution $u_f(t,.)$ and of the first term $H_f(t,.,u_0,t_0,x_0)$ for all $t \neq t_0$. We give firstly a formula for the difference between the two variances for arbitrary $(t_0,x_0)$ and we determine secondly the value of $(t_0,x_0)$ so that this difference is constant for all $t \neq t_0$.

\begin{4-PROP3} \label{4-PROP3}
	Let $(t_0,x_0) \in \R^2$. Suppose that the hypotheses of Theorem \ref{3-THM1} are satisfied and suppose in addition that $\| u_0 \|_{L^2(\R)} = 1$. Then for all $t \in \R \backslash \{ t_0 \}$, we have
	\begin{align}
		& \mathcal{V}\big( u_f(t,.) \big) - \mathcal{V}\big(H_f(t,.,u_0, t_0, x_0) \big) \nonumber \\
		& \hspace{1.5cm} =  2 \bigg( \frac{1}{2\pi} \, \Im \hspace{-1mm} \left( \int_{\R} f'(p) \, \tf u_0 (p) \, \overline{\big( \tf u_0 \big)'(p)} \, dp \right) - \mathcal{M}_{f'} \hspace{-1mm} \left( \frac{1}{\sqrt{2 \pi}} \tf u_0 \right) \mathcal{M}_1(u_0) \nonumber \\
		& \hspace{3cm} + \mathcal{V}_{f'}\hspace{-1mm} \left( \frac{1}{\sqrt{2 \pi}} \, \tf u_0 \right) t_0 \bigg) \, t + \mathcal{V}(u_0) - \mathcal{V}_{f'}\hspace{-1mm} \left( \frac{1}{\sqrt{2 \pi}} \, \tf u_0 \right) t_0^{\, 2} \; . \label{eq:var_diff}
	\end{align}
\end{4-PROP3}

\begin{proof}
	Let $t \in \R \backslash \{t_0\}$. Similarly to the arguments employed in the proof of Proposition \ref{4-PROP1}, we have
	\begin{align*}
		\mathcal{M}_2\big( H_f(t,.,u_0,t_0,x_0) \big)
			& = \frac{1}{2 \pi} \int_{x_0 + f'(p_1) (t-t_0)}^{x_0 + f'(p_2) (t-t_0)} x^2 \, \frac{ \big| \tf u_0 \big( p_0(t,x) \big) \big|^2}{f''\big(p_0(t,x) \big)} \, dx \, (t - t_0)^{-1} \\
			& = \frac{1}{2 \pi} \int_{p_1}^{p_2} \big(x_0 + f'(p) (t-t_0 ) \big)^2 \big| \tf u_0 (p) \big|^2 \, dp \\
			& = x_0^{\, 2} + \mathcal{M}_{f'^2} \hspace{-1mm} \left( \frac{1}{\sqrt{2 \pi}} \tf u_0 \right) (t-t_0)^2 + 2 \, x_0 \, \mathcal{M}_{f'} \hspace{-1mm} \left( \frac{1}{\sqrt{2 \pi}} \tf u_0 \right) (t-t_0) \; .
	\end{align*}
	Moreover, applying Proposition \ref{4-PROP1} gives
	\begin{equation*}
		\mathcal{M}_1\big( H_f(t,.,u_0,t_0,x_0) \big)^2 = x_0^{\, 2} + \mathcal{M}_{f'} \hspace{-1mm} \left( \frac{1}{\sqrt{2 \pi}} \, \tf u_0 \right)^2 (t-t_0)^2 + 2 \, x_0 \, \mathcal{M}_{f'} \hspace{-1mm} \left( \frac{1}{\sqrt{2 \pi}} \, \tf u_0 \right) (t-t_0) \; .
	\end{equation*}
	It follows:
	\begin{align*}
		\mathcal{V}\big(H_f(t,.,u_0,t_0,x_0) \big)
			& = \left( \mathcal{M}_{f'^2} \hspace{-1mm} \left( \frac{1}{\sqrt{2 \pi}} \, \tf u_0 \right) - \mathcal{M}_{f'} \hspace{-1mm} \left( \frac{1}{\sqrt{2 \pi}} \, \tf u_0 \right)^2 \right) (t-t_0)^2 \\
			& = \mathcal{V}_{f'}\hspace{-1mm} \left( \frac{1}{\sqrt{2 \pi}} \, \tf u_0 \right) (t-t_0)^2 \; .
	\end{align*}
	Using the formula for $\mathcal{V}\big( u_f(t,.) \big)$ from Proposition \ref{A-PROP2}, we obtain finally
	\begin{align*}
		& \mathcal{V}\big( u_f(t,.) \big) - \mathcal{V}\Big(H_f(t,.,u_0,t_0,x_0) \Big) \\
		& \hspace{1.5cm} = \mathcal{V}_{f'}\hspace{-1mm} \left( \frac{1}{\sqrt{2 \pi}} \, \tf u_0 \right) t^2 + 2 \bigg( \frac{1}{2\pi} \, \Im \hspace{-1mm} \left( \int_{\R} f'(p) \, \tf u_0 (p) \, \overline{\big( \tf u_0 \big)'(p)} \, dp \right) \\
		& \hspace{3cm} - \mathcal{M}_{f'} \hspace{-1mm} \left( \frac{1}{\sqrt{2 \pi}} \, \tf u_0 \right) \mathcal{M}_1(u_0) \bigg) t + \mathcal{V}(u_0) \\
		&\hspace{4.5cm} - \mathcal{V}_{f'}\hspace{-1mm} \left( \frac{1}{\sqrt{2 \pi}} \, \tf u_0 \right) (t-t_0)^2 \\
		& \hspace{1.5cm} =  2 \bigg( \frac{1}{2\pi} \, \Im \hspace{-1mm} \left( \int_{\R} f'(p) \, \tf u_0 (p) \, \overline{\big( \tf u_0 \big)'(p)} \, dp \right) - \mathcal{M}_{f'} \hspace{-1mm} \left( \frac{1}{\sqrt{2 \pi}} \tf u_0 \right) \mathcal{M}_1(u_0) \\
		& \hspace{3cm} + \mathcal{V}_{f'}\hspace{-1mm} \left( \frac{1}{\sqrt{2 \pi}} \, \tf u_0 \right) t_0 \bigg) \, t + \mathcal{V}(u_0) - \mathcal{V}_{f'}\hspace{-1mm} \left( \frac{1}{\sqrt{2 \pi}} \, \tf u_0 \right) t_0^{\, 2} \; .
	\end{align*}
\end{proof}

In view of the preceding result, the difference between the variances of $u_f(t,.)$ and $H_f(t,.,u_0,t_0,x_0)$ is an affine function with respect to $t$. Consequently the unique way to make this difference constant is to choose $t_0$ in such way that the leading coefficient is equal to $0$. It turns out that the unique $t_0$ satisfying this property is the one minimizing the variance of $u_f$, namely $t^*$ introduced in Lemma \ref{3-LEM1}.

\begin{4-PROP4} \label{4-PROP4}
	Let $(t_0,x_0) \in \R^2$. Suppose that the hypotheses of Theorem \ref{3-THM1} are sa\-tisfied and suppose in addition that $\| u_0 \|_{L^2(\R)} = 1$. Then we have the following equivalence:
	\begin{align*}
		& \exists \, C \in \R \quad \forall \, t \in \R \backslash \{t_0\}  \qquad \mathcal{V} \big( u_f (t,.) \big) - \mathcal{V}\big( H_f(t,.,u_0,t_0,x_0) \big) = C \\
		& \hspace{1.5cm} \Longleftrightarrow \qquad t_0 = t^* := \argmin_{\tau \in \R} \mathcal{V}\big( u_f(\tau,.) \big) \; .
	\end{align*}
	In particular, we have
	\begin{equation*}
		C = \min_{\tau \in \R} \mathcal{V}\big(u_f(\tau,.) \big) \; .
	\end{equation*}
\end{4-PROP4}

\begin{proof}
	According to Proposition \ref{4-PROP3}, the difference between the variances of $u_f(t,.)$ and $H_f(t,.,u_0,t_0,x_0)$ is constant if and only if
	\begin{align*}
		t_0
			& = \frac{1}{\mathcal{V}_{f'}\big( \frac{1}{\sqrt{2 \pi}} \, \tf u_0 \big)} \left( - \frac{1}{2 \pi} \Im \bigg(  \int_{\R} f'(p) \, \tf u_0(p) \, \overline{\big( \tf u_0 \big)'(p)} \, dp \right) \\
			& \hspace{1.5cm} + \mathcal{M}_{f'} \hspace{-1mm} \left( \frac{1}{\sqrt{2 \pi} } \, \tf u_0 \right) \mathcal{M}_1(u_0) \bigg) \; ,
	\end{align*}
	which is equal to $t^*$ according to Remark \ref{3-REM1}. By evaluating equality \eqref{eq:var_diff} at $t^*$, we obtain for all $t \in \R \backslash \{t^* \}$,
	\begin{equation} \label{eq:a}
		\mathcal{V}\big( u_f(t,.) \big) - \mathcal{V}\big(H_f(t,.,u_0,t^*, x_0 \big) = \mathcal{V}(u_0) - \mathcal{V}_{f'}\hspace{-1mm} \left( \frac{1}{\sqrt{2 \pi}} \, \tf u_0 \right) (t^*)^2 \; .
	\end{equation}
	We note now that
	\begin{align}
		\mathcal{V}\big(u_f(t^*,.) \big)
			& = \mathcal{V}_{f'}\hspace{-1mm} \left( \frac{1}{\sqrt{2 \pi}} \, \tf u_0 \right) (t^*)^2 + 2 \bigg( \frac{1}{2\pi} \, \Im \hspace{-1mm} \left( \int_{\R} f'(p) \, \tf u_0 (p) \, \overline{\big( \tf u_0 \big)'(p)} \, dp \right) \nonumber \\
		& \hspace{3cm} - \mathcal{M}_{f'} \hspace{-1mm} \left( \frac{1}{\sqrt{2 \pi}} \, \tf u_0 \right) \mathcal{M}_1(u_0) \bigg) t^* + \mathcal{V}(u_0) \nonumber \\
			& = \mathcal{V}_{f'}\hspace{-1mm} \left( \frac{1}{\sqrt{2 \pi}} \, \tf u_0 \right) (t^*)^2 - 2 \mathcal{V}_{f'}\hspace{-1mm} \left( \frac{1}{\sqrt{2 \pi}} \, \tf u_0 \right) (t^*)^2 + \mathcal{V}(u_0) \nonumber \\
			& = - \mathcal{V}_{f'}\hspace{-1mm} \left( \frac{1}{\sqrt{2 \pi}} \, \tf u_0 \right) (t^*)^2 + \mathcal{V}(u_0) \; . \label{eq:b}
	\end{align}
	From equalities \eqref{eq:a} and \eqref{eq:b}, it follows finally
	\begin{equation*}
		\mathcal{V}\big( u_f(t,.) \big) - \mathcal{V}\big(H_f(t,.,u_0,t^*, x_0 \big) = \mathcal{V}\big(u_f(t^*,.) \big) = \min_{\tau \in \R} \mathcal{V}\big(u_f(\tau,.) \big) \; ,
	\end{equation*}
	the last equality being obtained by the definition $\displaystyle t^* = \argmin_{\tau \in \R} \mathcal{V}\big(u_f(\tau,.) \big)$.
\end{proof}

According to Propositions \ref{4-PROP2} and \ref{4-PROP4}, choosing $(t^*, x^*)$, introduced in Lemma \ref{3-LEM1}, as the origin of the cone in which we expand the solution $u_f(t,.)$ of equation \eqref{eq:evoleq} provides a time-asymptotic approximation $H_f(t,.,u_0,t^*,x^*)$ having the right mean position and a constant error. This is summarized in the following corollary.

\begin{4-COR1}
	Suppose that the hypotheses of Theorem \ref{3-THM1} are satisfied and suppose in addition that $\| u_0 \|_{L^2(\R)} = 1$. Then for all $t \in \R \backslash \{ t^* \}$, we have
	\begin{equation*}
		\left\{ \begin{array}{l}
			\displaystyle \mathcal{M}_1 \big( u_f (t,.) \big) = \mathcal{M}_1\big( H_f(t,.,u_0,t^*, x^*) \big) \\[2mm]
			\displaystyle \mathcal{V} \big( u_f (t,.) \big) - \mathcal{V}\big( H_f(t,.,u_0,t^*,x^*) \big) = \min_{\tau \in \R} \mathcal{V}\big(u_f(\tau,.) \big)
		\end{array} \right. \; ,
	\end{equation*}
	where $t^*$ and $x^*$ are defined in Lemma \ref{3-LEM1}.
\end{4-COR1}

\begin{proof}
	Simple application of Propositions \ref{4-PROP2} and \ref{4-PROP4} combined with the definitions of $t^*$ and $x^*$.
\end{proof}

\section{Appendix A: Mean position and variance of the free wave packet}

In this appendix, we give the formulas for the mean position and the variance of the wave packet $u_f$ defined in \eqref{eq:sol_formula0}. The proofs we propose here are substantially based on the fact that the wave packet is defined via the Fourier transform, permitting to apply some pro\-per\-ties of the Fourier transform.\\

We begin with the formula for the mean position.

\begin{4-PROP1} \label{A-PROP1}
	Suppose that $u_0 \in \mathcal{S}(\R)$. Then for all $t \in \R$, we have
	\begin{equation*}
		\mathcal{M}_1 \big( u_f(t,.) \big) = \mathcal{M}_{f'} \hspace{-1mm} \left( \frac{1}{\sqrt{2 \pi}} \, \tf u_0 \right) t + \mathcal{M}_1(u_0) \; .
	\end{equation*}
\end{4-PROP1}

\begin{proof}
	For $t \in \R$, we have
	\begin{align}
		\int_{\R} x \, \big| u_f(t,x) \big|^2 \, dx & = \int_{\R} x \, u_f(t,x) \, \overline{u_f(t,x)} \, dx \nonumber \\
			& = \frac{1}{2 \pi} \int_{\R} \tf\big[x \mapsto x \, u_f(t,x)\big](p) \, \overline{\tf \big[x \mapsto u_f(t,x) \big](p)} \, dp \nonumber \\
			& = \frac{i}{2 \pi} \int_{\R} \partial_p \tf\big[x \mapsto u_f(t,x)\big](p) \, \overline{\tf \big[x \mapsto u_f(t,x) \big](p)} \, dp \; ; \label{eq:ehren}
	\end{align}
	the second and third equalities have been obtained by applying Plancherel theorem and basic properties of the Fourier transform. Using now the expression \eqref{eq:sol_formula0} of the wave packet $u_f(t,x)$, we obtain for all $p \in \R$,
	\begin{align*}
		& \bullet \quad \partial_p \tf\big[x \mapsto u_f(t,x)\big](p) = e^{-itf(p)} \Big( -i t f'(p) \, \tf u_0(p) + (\tf u_0)'(p) \Big) \; ; \\[2mm]
		& \bullet \quad \overline{\tf \big[x \mapsto u_f(t,x) \big](p)} = e^{itf(p)} \, \overline{\tf u_0(p)} \; .
	\end{align*}
	By combing the two last equalities with \eqref{eq:ehren} and by using again basic properties of the Fourier transform, it follows
	\begin{align*}
		\int_{\R} x \, \big| u_f(t,x) \big|^2 \, dx & = \frac{i}{2 \pi} \int_{\R} \Big( -i t f'(p) \, \tf u_0(p) \: + \: (\tf u_0)'(p) \Big) \, \overline{\tf u_0(p)} \, dp \\
			& = \frac{1}{2 \pi} \int_{\R} f'(p) \, \big| \tf u_0 (p) \big|^2 \, dp \, t \: + \: \frac{i}{2 \pi} \int_{\R} (\tf u_0)'(p) \, \overline{\tf u_0(p)} \, dp \\
			& = \frac{1}{2 \pi} \int_{\R} f'(p) \, \big| \tf u_0 (p) \big|^2 \, dp \, t \: + \: \frac{1}{2 \pi} \int_{\R} \tf \big[x \mapsto x \, u_0(x)\big](p) \, \overline{\tf u_0(p)} \, dp \\
			& = \frac{1}{2 \pi} \int_{\R} f'(p) \, \big| \tf u_0 (p) \big|^2 \, dp \, t \: + \: \int_{\R} x \, \big| u_0(x) \big|^2 \, dx \; ,
	\end{align*}
	leading finally to the desired equality.
\end{proof}

\begin{A-REM1} \em
	The preceding formula is actually an extension of the well-known Ehrenfest theorem to the family of dispersive equations of type \eqref{eq:evoleq0}. We recall that Ehrenfest theorem in the setting of the free Schrödinger equation \eqref{eq:free_schrodinger} gives the following formula for the mean position of the free particle:
	\begin{equation*}
		\mathcal{M}_1\big( u_{S}(t,.) \big) = \mathcal{M}_1 \hspace{-1mm} \left( \frac{1}{\sqrt{2 \pi}} \, \tf u_0 \right) t + \mathcal{M}_1(u_0) \; .
\end{equation*}
\end{A-REM1}

The formula for the variance is provided in the following result.

\begin{A-PROP2} \label{A-PROP2}
	Suppose that $u_0 \in \mathcal{S}(\R)$. Then for all $t \in \R$, we have
	\begin{align*}
		\mathcal{V}\big( u_f(t,.) \big)
			& = \mathcal{V}_{f'}\hspace{-1mm} \left( \frac{1}{\sqrt{2 \pi}} \, \tf u_0 \right) t^2 + 2 \bigg( \frac{1}{2\pi} \, \Im \hspace{-1mm} \left( \int_{\R} f'(p) \, \tf u_0 (p) \, \overline{\big( \tf u_0 \big)'(p)} \, dp \right) \\
			& \hspace{1.5cm} - \mathcal{M}_{f'} \hspace{-1mm} \left( \frac{1}{\sqrt{2 \pi}} \, \tf u_0 \right) \mathcal{M}_1(u_0) \bigg) t + \mathcal{V}(u_0) \; .
	\end{align*}
\end{A-PROP2}

\begin{proof}
	Following the computational arguments of the proof of Proposition \ref{A-PROP1}, we have for all $t \in \R$,
	\begin{align}
		\int_{\R} x^2 \, \big| u_f(t,x) \big|^2 \, dx & = \int_{\R} \big| x \, u_f(t,x) \big|^2 \, dx \nonumber \\
			& = \frac{1}{2 \pi} \int_{\R} \Big| \tf\big[x \mapsto x \, u_f(t,x)\big](p) \Big|^2 \, dp \nonumber \\
			& = \frac{1}{2 \pi} \int_{\R} \Big| \partial_p \tf\big[x \mapsto u_f(t,x)\big](p) \Big|^2 \, dp \nonumber \\
			& =  \frac{1}{2 \pi} \int_{\R} \Big| - i t f'(p) \, \tf u_0(p) + (\tf u_0)'(p) \Big|^2 \, dp \; , \label{eq:var}
	\end{align}
	and we recall that
	\begin{equation*}
		\forall \, p \in \R \qquad (\tf u_0)'(p) = -i \, \tf \big[ x \mapsto x \, u_0(x) \big](p) \; .
	\end{equation*}
	Inserting the preceding relation into \eqref{eq:var} and expanding then the square of the absolute value provides
	\begin{align*}
		\int_{\R} x^2 \, \big| u_f(t,x) \big|^2 \, dx
			& = \frac{1}{2\pi} \int_\R f'(p)^2 \, \big| \tf u_0 (p) \big|^2 \, dp \, t^2 \: + \: \int_\R x^2 \, \big| u_0 (x) \big|^2 \, dx \nonumber \\
			& \hspace{1cm} - \: \frac{1}{\pi} \int_{\R} \Re \hspace{-1mm} \left( i \, f'(p) \, \tf u_0 (p) \, \overline{\big( \tf u_0 \big)'(p)} \right) dp \, t \\
			& = \mathcal{M}_{f'^2} \hspace{-1mm} \left( \frac{1}{\sqrt{2 \pi}} \, \tf u_0 \right) t^2 \: + \: \mathcal{M}_2(u_0) \nonumber \\
			& \hspace{1cm} + \: \frac{1}{\pi} \, \Im \hspace{-1mm} \left( \int_{\R} f'(p) \, \tf u_0 (p) \, \overline{\big( \tf u_0 \big)'(p)} \, dp \right) t \; ,
	\end{align*}
	Now by using Proposition \ref{A-PROP1}, we have
	\begin{equation*}
		\mathcal{M}_1 \big( u_f(t,.) \big)^2 = \mathcal{M}_{f'} \hspace{-1mm} \left( \frac{1}{\sqrt{2 \pi}} \, \tf u_0 \right)^{\hspace{-1mm} 2} t^2 + \mathcal{M}_1(u_0)^2 + 2 \, \mathcal{M}_{f'} \hspace{-1mm} \left( \frac{1}{\sqrt{2 \pi}} \, \tf u_0 \right) \mathcal{M}_1(u_0) \, t \; .
	\end{equation*}
	which leads finally to
	\begin{align*}
		\mathcal{V}\big( u_f(t,.) \big)
			&= \mathcal{M}_2\big( u_f(t,.) \big) - \mathcal{M}_1\big( u_f(t,.) \big)^2 \\
			& = \mathcal{M}_{f'^2} \hspace{-1mm} \left( \frac{1}{\sqrt{2 \pi}} \, \tf u_0 \right) t^2 \: + \: \mathcal{M}_2(u_0) \: + \: \frac{1}{\pi} \, \Im \hspace{-1mm} \left( \int_{\R} f'(p) \, \tf u_0 (p) \, \overline{\big( \tf u_0 \big)'(p)} \, dp \right) t \\
		& \hspace{1cm} - \mathcal{M}_{f'} \hspace{-1mm} \left( \frac{1}{\sqrt{2 \pi}} \, \tf u_0 \right)^{\hspace{-1mm} 2} t^2 - \mathcal{M}_1(u_0)^2 \: - \: 2 \, \mathcal{M}_{f'} \hspace{-1mm} \left( \frac{1}{\sqrt{2 \pi}} \, \tf u_0 \right) \mathcal{M}_1(u_0) \, t \\
			& = \mathcal{V}_{f'}\hspace{-1mm} \left( \frac{1}{\sqrt{2 \pi}} \, \tf u_0 \right) t^2 + \mathcal{V}(u_0) + 2 \bigg( \frac{1}{2\pi} \, \Im \hspace{-1mm} \left( \int_{\R} f'(p) \, \tf u_0 (p) \, \overline{\big( \tf u_0 \big)'(p)} \, dp \right) \\
			& \hspace{1cm} - \mathcal{M}_{f'} \hspace{-1mm} \left( \frac{1}{\sqrt{2 \pi}} \, \tf u_0 \right) \mathcal{M}_1(u_0) \bigg) t \; .
	\end{align*}
\end{proof}

\noindent \textbf{Acknowledgements:}\\
The author gratefully thanks Prof.~Felix Ali Mehmeti whose numerous and insightful comments helped to improve the present paper.

\end{document}